\documentclass[article,12pt]{elsarticle}

\usepackage{amssymb}
\usepackage{bm}
\usepackage{mathrsfs}
\usepackage{multirow}
\usepackage{chngpage}
\usepackage{array}
\usepackage{rotating}
\usepackage{setspace}
\usepackage{color}
\usepackage{geometry}
\usepackage{mathrsfs}
\usepackage{amsmath}
\begin{document}
\doublespacing

\begin{frontmatter}



\title{Semiparametric model averaging for high dimensional conditional quantile prediction}

\author{Jingwen Tu$^{a}$, Hu Yang$^{a}$, Chaohui Guo$^{b*}$}

\address{
$^a$College of Mathematics and Statistics, Chongqing University, Chongqing, 401331, China
\\
$^b$College of Mathematics Science, Chongqing Normal University, Chongqing, 401331, China }
\cortext[*]{Corresponding author (Email:guochaohui2010@126.com)}

\begin{abstract}
In this article, we propose a penalized high dimensional semiparametric model average quantile prediction approach that is robust for forecasting the conditional quantile of the response. We consider a two-step estimation procedure. In the first step, we use a local linear regression approach to estimate the individual marginal quantile functions, and approximate the conditional quantile of the response by an affine combination of one-dimensional marginal quantile regression functions. In the second step, based on the nonparametric kernel estimates of the marginal quantile regression functions, we utilize a penalized method to estimate the suitable model weights vector involved in the approximation. The objective of the second step is to select significant variables whose marginal quantile functions make a significant contribution to estimating the joint multivariate conditional quantile function. Under some mild conditions, we have established the asymptotic properties of the proposed robust estimator. Finally, simulations and a real data analysis have been used to illustrate the proposed method.

\end{abstract}


\end{frontmatter}

% \linenumbers


\section{Introduction}

In many practical situations, especially for economic and medical fields, forecasting and predictive inference are our main goals. In practice, we often face a large number of predictors and uncertain functional forms when making statistical prediction. A popular approach to solve this problem is to consider the model selection tool that can select a optimal model from all candidate models, but we have to recognize that model selection technique yields only one final model, so useful information may be ignored when significant variables absent from the final model. This may result in misleading predictive outcomes. Instead of depending on only one best model, an alternative method,  called model averaging technique, aims to improve the prediction accuracy through giving higher weights to the better marginal models.  Thus, model averaging can be regarded as a smoothed extension of model selection and generally leads to a lower risk than model selection. Earlier development for model average was linked closely the Bayesian statistics including Hoeting et al. (1999), Raftery et al. (1997) and Hjort and Claeskens (2003). Recently, various strategies have been developed to construct optimal model averaging weights for frequentist models. For example, Hansen (2007) proposed a frequentist model average approach with weights selected by minimizing a Mallows criterion. Wan et al. (2010) focused on two assumptions of Hansen (2007) and provided a stronger theoretical basis for the use of the Mallows criterion in model averaging. Liang et al. (2011) considered a new procedure of weight choice by minimizing frequentist model average estimators' mean squared errors. To deal with heteroscedastic data, Hansen and Racine (2012) developed a jackknife model averaging approach to choose weights by minimizing a leave-one-out cross-validation criterion and had proved that the proposed approach achieved the lowest possible asymptotic squared error. Zhang et al. (2013) further extended the method of Hansen and Racine (2012) to general models with a non-diagonal error covariance structure or lagged dependent variables. In the framework of linear mixed-effects models, Zhang et al. (2014) constructed an unbiased estimator of the squared risk for the model averaging, which has been demonstrated to be asymptotically optimal in theory under some regularity conditions. Zhang et al. (2016) studied optimal model averaging methods for generalized linear models and generalized linear mixed-effects models, which can be taken as an extension of Zhang et al. (2014)'s. Under the local asymptotic framework, Liu et al. (2015) studied the limiting distributions of least squares averaging estimators and proposed a plug-in averaging estimator by minimizing the sample asymptotic mean squared error. Other related literature can refer to Hansen (2008), Claeskens and Hjort (2008), Zhang et al. (2012), Cheng and Hansen (2015).

Almost all mentioned above research work focus on averaging a set of parametric models by assuming some parametrically linear or nonlinear relationships between the response and predictors. Although parametric models are easy to understand and widely accepted
by scientific researchers, they make strong assumptions in practical applications, which may increase the risk of bias prediction.
In contrast, nonparametric models with less structural restriction may provide more flexible predictive inference. Recently, Li et al. (2015) firstly proposed a nonparametric model averaging approach which is more flexible than traditional parametric averaging method. They estimated the multivariate conditional mean regression function by averaging a set of estimated marginal mean regression functions with proper weights obtained by minimizing least squares loss. Motivated by the nonparametric model averaging technique, Chen et al. (2016) studied the semiparametric dynamic portfolio choice and utilized a novel data-driven method to estimate the nonparametric optimal portfolio choice. Huang and Li (2018) extended the method of Li et al. (2015) to panel data and established the asymptotic results of the proposed procedure. Li et al. (2018) approximated the conditional mean regression function by a weighted average of varying coefficient regression functions, which can handle discrete and continuous predictors.

In recent years, we often encounter datasets with a very large number of potential predictors, but only a minority of predictors are truly relevant in prediction. However, most of literature focus on the determination of weights for individual models under a fixed number of covariates. So far, Ando and Li (2014) proposed a two-step model averaging procedure to predict the conditional mean of the response for a
ultra-high dimensional linear regression. In order to obtain more accurate prediction of the conditional mean of the response for ultra-high dimensional time series, Chen et al. (2018) introduced a two-step semiparametric procedure that includes the kernel sure independence screening technique and the semiparametric penalized method of model averaging marginal regression. All mentioned above references aim to forecast the conditional mean of the response, but sometimes we are more interested in predicting the conditional quantile of the response. Compared to mean regression, quantile regression not only provides a more complete description of
the entire response distribution but also does not require specification of the error distribution, and thus it is more robust.

In this paper, we aims to develop a new semiparametric model averaging procedure for achieving more accurate prediction for the true conditional quantile of the response under the high dimensional setting. This paper may have several innovation as follows: (1)
our objective is to predict the conditional quantile of the response rather than its conditional mean. Thus we may encounter more challenge to establish asymptotic theories of model weights since we cannot obtain the closed-form expression of model weights; (2) the proposed approach can offer a complete prediction for the response when different quantiles are adopted; (3) our method produces more accurate in-sample and out-of-sample prediction when non-normal error are considered.

The rest of the paper is organized as follows. In Section 2, we first give the approximation of the conditional quantile function of the response. Then a two-step semiparametric model averaging approach is applied to estimate the conditional quantile function of the response.
In Section 3, we establish the asymptotic theory for the proposed estimator. In Section 4, numerical studies including simulation studies and a real data analysis are carried out to investigate the finite sample performance of the proposed method. Some discussions are reported in Section 5. Finally, all technical proofs are given in the Appendix.

\section{Model approximation and estimation method}

Let $\{ (\bm{X}_i,Y_i),1 \le i \le n\} $ be independent and identically distributed observations from $\left(\bm X,Y\right ) $, where $\bm{X}=\left( X_1,...,X_{p_n}  \right)$ is a $p_n$-vector of predictors and $Y$ is the response variable. The goal of this paper is to develop new procedure for forecasting the $\tau$th conditional quantile function of $Y$ given $\bm X$, namely, $m_{\tau}\left(\bm X\right) \equiv Q_{\tau}\left( {Y|\bm X} \right)$. If the dimension of $\bm X$ is high, it is not practical to model conditional quantile function $ m_{\tau}\left(\bm X\right)$ without any structure assumption due to the curse of dimensionality. Recently, authors approximated the quantile function $m_{\tau}\left(\bm X\right)$ by semiparametric models such as quantile additive models (Horowitz and Lee, 2005, Lv et al., 2017), quantile varying coefficient models (Tang et al., 2013) and among others. However, using a specified model with fixed model structure may increase the risk of model misspecification, which results in poor predictive performance. Therefore, we adopt the model averaging technique to predict $m_{\tau}\left(\bm X\right)$. Specifically, motivated by Li et al. (2015), we model or approximate $m_{\tau}\left(\bm X\right)$ by an affine combination of one-dimensional nonparametric functions ${m_\tau }\left( \bm X \right) = {w_0} + \sum\nolimits_{j = 1}^{p_n} {w_j{m_{\tau_j j}}\left( {{X_j}} \right)}$, where $m_{\tau_jj}\left( {{X_j}} \right)={Q_{\tau_j} }\left( {Y|{X_j}} \right)$
is the $\tau$th conditional quantile of $Y$ given $X_j$. Here, each marginal regression $m_{\tau j}\left( \cdot \right)$ can be regarded as a candidate model and $w_j$ is the corresponding model weight coefficient. In the rest of the article, we omit $\tau_j$ and $\tau$ from $m_{\tau_j j}(\cdot)$ and $m_{\tau}\left(\bm X\right)$ for notational simplicity, but it is helpful to bear in mind that these quantities are $\tau_j$ and $\tau$-specific.

What we are most interested in is to accurately estimate $m_{j} $ and the model average weight vector $\bm{{\rm{ w}}}=\left( w_0,w_1,...,w_{p_n}\right)^T$. We consider a two-step estimation procedure. In the first step, we employ local linear regression technique to estimate the individual marginal regression functions $m_{j}\left( \cdot \right), j=1,...,p_n$. Specifically, considering a Taylor
expansion, we have
\[m_{ j}\left( {{X_{ij}}} \right) \approx {m_{ j}}\left( x \right) + {{\dot m}_{ j}}\left( x \right)\left( {{X_{ij}} - x} \right) \equiv a + b\left( {{X_{ij}} - x} \right),i=1,...,n, j=1,...,p_n,\]
where ${{\dot m}_{j}}\left( x \right)$ is the first-order derivative of ${{ m}_{ j}}\left( x \right)$. Let ${\rho _\tau }\left( u \right) = \tau u - uI\left( {u < 0} \right)$ be check loss function at $\tau $ quantile. Then, we estimate $m_{j}$ by minimizing the following local weighted quantile loss
\begin{equation}
\sum\limits_{i = 1}^n {{\rho _{\tau_j} }\left\{ {{Y_i} - a - b\left( {{X_{ij}} - x} \right)} \right\}} K\left( {\frac{{{X_{ij}} - x}}{h_j}} \right),
\end{equation}
where $K(\cdot)$ is a kernel function and $h_j$ is a bandwidth. Let $\left( {\hat a,\hat b} \right)$ be the minimizer of the objective function (1). Then, we have ${{\hat m}_{j}} = \hat a$.

In the second step, let ${\bm{{\rm{w}}}_o}  = {\left( {{w_{o0}},{w_{o1}},...,{w_{o{p_n}}}} \right)^T}$ be the optimal values of the weights in the model averaging defined in Li et al. (2015). To estimate ${\bm{{\rm{w}}}_o}$, we minimize the following function with respect to $ {\bm{{\rm{w}}}_{n}}  = \left( {{w_{0}},\bm w^T} \right)^T$ with $\bm w={\left( {{w_{1}},...,{w_{{p_n}}}} \right)^T}$,
\begin{equation}
\mathcal{Q}_n\left( {\bm{{\rm{w}}}_{n}}\right)= \sum\limits_{i = 1}^n {{\rho _\tau }\left\{ {{Y_i} - w_0-\sum\limits_{j = 1}^{p_n} {{{\hat m}_j}\left( {{X_{ij}}} \right){w _j}} } \right\}} +n\sum\limits_{j = 1}^{{p_n}} {{p_\lambda }\left( {\left| {{w_j}} \right|} \right)} ,
\end{equation}
where ${p_\lambda }\left(\cdot\right)$ is a penalty function with a tuning parameter $\lambda$, such as SCAD penalty function, $\dot p _{\lambda}(\cdot)$ is its first order derivative, defined by
\[{\dot{p}_{\lambda}}(x) = {\lambda}\left\{ {I\left( {x \le \lambda} \right) + \frac{{{{\left( {a\lambda - x} \right)}_ + }}}{{\left( {a - 1} \right)\lambda}}I\left( {x > \lambda} \right)} \right\},\]
where $a > 2$, ${{p}_{\lambda}}(0)=0$ and $\lambda$ is a nonnegative penalty parameter which governs sparsity of the model. It is easy to find that $\dot{p}_\lambda(\left| x \right|) $ is close to zero if $\left| x \right|$ is large.

The estimator of the optimal weights ${\bm{{\rm{w}}}_o}$ can be obtained through minimizing the objective function (2), that is, $ {\bm{{\rm{\hat w}}}_{n}}= \mathop {\arg \min }\limits_{ {\bm{{\rm{w}}}_{n}}}{\mathcal{Q}_n\left( {\bm{{\rm{w}}}_{n}}\right)} $. This paper uses the R package ``rqPen'' to obtain the estimator $ {\bm{{\rm{\hat w}}}_{n}}$. Finally, for a future observation $\bm{{\rm{x}}}=\left({\rm{x}}_1,...,{\rm{x}}_{p_n}\right)$, we can predict $m(\bm{{\rm{x}}})$ by $\hat {m}\left( \bm{{\rm{x}}}  \right) = \hat{w}_0 + \sum\nolimits_{j = 1}^{p_n} {\hat w_j{\hat m_{j}}\left( {\rm{x}}_j \right)}$.

\section{The theoretical results}

Define ${a_n} = \mathop {\max }\limits_{1 \le j \le {p_n}} \left\{ {\left| {{{\dot p}_\lambda }\left( {\left| {{w_{oj}}} \right|} \right)} \right|,\left| {{w_{oj}}} \right| \ne 0} \right\}$ and ${b_n} = \mathop {\max }\limits_{1 \le j \le {p_n}} \left\{ {\left| {{{\ddot p}_\lambda }\left( {\left| {{w_{oj}}} \right|} \right)} \right|,\left| {{w_{oj}}} \right| \ne 0} \right\}$, where ${{\ddot p}_{\lambda}}\left( x \right)$ is the second-order derivative of ${{ p}_{ \lambda}}\left( x \right)$, ${\eta _{ij}} = {Y_i} - {m_j}\left( {{X_{ij}}} \right)$, ${\eta _i} = {Y_i} - {w_{o0}} - \sum\limits_{j = 1}^{{p_n}} {{w_{oj}}{m_j}\left( {{X_{ij}}} \right)}$, $\mathcal{M}_i = {\left( {{m_1}\left( {{X_{i1}}} \right),...,{m_{{p_n}}}\left( {{X_{i{p_n}}}} \right)} \right)^T}$ and ${{\hat {\mathcal{M}}}_i} = {\left( {{{\hat m}_1}\left( {{X_{i1}}} \right),...,{{\hat m}_{p_n}}\left( {{X_{i{p_n}}}} \right)} \right)^T}$. We assume $p\left\{\eta_{ij}\le 0 \right\}=\tau_j$ and $p\left\{\eta_{i}\le 0 \right\}=\tau$. To prove the theoretical results of the proposed estimators, we next present the following technical conditions.

\noindent (C1) Let $f_j(\cdot)$ be the marginal density function of the covariates $\left\{ {{X_{ij}}} \right\}$, the $j$-th element of $\bm X_i$. Assume that $f_j(\cdot)$ has continuous derivatives up to the second order and
\[0 <  c\le \mathop {\inf }\limits_j \mathop {\inf }\limits_{{x_j} \in {\mathcal{C}_j}} {f_j}\left( {{x_j}} \right) \le \mathop {\sup }\limits_j \mathop {\sup }\limits_{{x_j} \in {\mathcal{C}_j}} {f_j}\left( {{x_j}} \right) \le C< \infty \]
where $\mathcal{C}_j$ is the compact support of $X_{ij}$. For each $j$, the conditional density functions of $Y_i$ for
given $X_{ij}$ exists and satisfies the Lipschitz continuous condition. Furthermore, the length of $\mathcal{C}_j$
is uniformly bounded by a positive constant.

\noindent (C2) The kernel function $K\left(  \cdot  \right)$ is a Lipschitz continuous, symmetric and bounded probability density
function with a compact support.

\noindent (C3) The marginal regression function $m_j\left(  \cdot  \right)$ has continuous derivatives up to the second order and
there exists a positive constant $c_m$ such that
 \[{\sup _j}{\sup _{{x_j} \in {\mathcal{C}_j}}}\left[ {\left| {{m_j}\left( {{x_j}} \right)} \right| + \left| {{{\dot m}_j}\left( {{x_j}} \right)} \right| + \left| {{{\ddot m}_j}\left( {{x_j}} \right)} \right|} \right] \le c_m.\]

\noindent (C4) Let $f_{\eta_j}(\cdot)$ and $F_{\eta_j}(\cdot)$ be the marginal density and distribution functions of $\eta_{ij}$, $f_{\eta}(\cdot)$ and $F_{\eta}(\cdot)$ be the density and distribution functions of $\eta_{i}$. The density functions $f_{\eta_j}(\cdot)$ and $f_{\eta}(\cdot)$ are bounded and bounded away from zero in a neighborhood of zero.

 \noindent (C5) There exists a sequence of fixed vectors $\left\{ \bm u \right\}$ in $\mathcal{R}^{p_n}$, with $\left\| \bm u \right\|$ bounded, such
that ${\max _{1 \le i \le n}}\left| {{\hat {\mathcal{M}}}_i^T\bm u} \right| = O\left( {\sqrt {\log n} } \right)$, where $\left\| \cdot \right\|$ denotes the $L_2$ norm for any vector.

 \noindent (C6) The matrix
\[{\bm\Lambda _n} = \left( {\begin{array}{*{20}{c}}
   {E\left[ {{m_1}\left( {{X_{i1}}} \right){m_1}\left( {{X_{i1}}} \right)} \right]} &  \cdots  & {E\left[ {{m_1}\left( {{X_{i1}}} \right){m_{{p_n}}}\left( {{X_{i{p_n}}}} \right)} \right]}  \\
    \vdots  &  \vdots  &  \vdots   \\
   {E\left[ {{m_{{p_n}}}\left( {{X_{i{p_n}}}} \right){m_1}\left( {{X_{i1}}} \right)} \right]} &  \cdots  & {E\left[ {{m_{{p_n}}}\left( {{X_{i{p_n}}}} \right){m_{{p_n}}}\left( {{X_{i{p_n}}}} \right)} \right]}  \\
\end{array}} \right)\]
is positive definite with the eigenvalues bounded away from zero and infinity. In particular, the
smallest eigenvalue of ${\bm\Lambda _n}$ is larger than $b$, a small positive constant.

\noindent (C7) $Var\left( {{{\hat {\mathcal{M}}}_i}} \right) = {\bm G_n} > 0,E\left( {\hat{ {\mathcal{M}}}_i^{ \otimes 2}} \right) = {\bm\Psi _n}$, and $0 < {d_1} \le {\lambda _{\min }}\left( {{\bm\Psi _n}} \right) \le {\lambda _{\max }}\left( {{\bm\Psi _n}} \right) < {d_2} < \infty $ for all $n$, where ${\lambda _{\min }}\left( {{\bm\Psi _n}} \right) $ and $ {\lambda _{\max }}\left( {{\bm\Psi _n}} \right)$ are the smallest and largest eigenvalues of ${\bm\Psi _n}$.

\noindent (C8) $\lim {\inf _{n \to \infty }}\lim {\inf _{x \to {0_ + }}}{\dot{p}_{\lambda}}\left( x \right)/\lambda > 0$.

\noindent (C9) Let $a_n=O(n^{-1/2})$ and $b_n=o(1)$, and there exist two positive constants $C_1$ and $C_2$ such that $\left| {{\ddot{p}_{\lambda}}\left( {{x_1}} \right) - {\ddot{p}_{\lambda}}\left( {{x_2}} \right)} \right| \le {C_2}\left| {{x_1} - {x_2}} \right|$ when $x_1,x_2>C_1\lambda$.

Without loss of generality, we define the vector of the optimal weights
\[{\bm{{\rm{w}}}_o}  = {\left( {{w_{o0}},{w_{o1}},...,{w_{o{p_n}}}} \right)^T}= \left({w_{o0}}, {\bm{{\rm{w}}}_o^T(1)}, {\bm{{\rm{w}}}_o^T(2)}    \right)^T,\]
where $ {\bm{{\rm{w}}}_o(1)} =\left( {{w_{o1}},...,{w_{o{s_n}}}} \right)^T$ stands for non-zero weights with dimension $s_n$ and $ {\bm{{\rm{w}}}_o(2)}=\left( {{w_{o,{s_n+1}}},...,{w_{o{p_n}}}} \right)^T $ is zero weights
with dimension $p_n-s_n$. Let ${\bm{{\rm{\hat w}}}_{n}(1)}$ and ${\bm{{\rm{\hat w}}}_{n}(2)}$ be the estimators of ${\bm{{\rm{ w}}}_{o}(1)}$ and ${\bm{{\rm{ w}}}_{o}(2)}$ respectively.

Define ${{\bar \eta }_{ij}} = {m_j}\left( {{X_{ij}}} \right){\psi _\tau }\left( {{\eta _i}} \right)$, ${\bm \mu _n} = E\left( {{{\hat {\mathcal{M}}}_{1i}}}\right)$, ${\bm\mu _{n1}}$ is first $s_n\times 1$ vector of ${\bm\mu _{n}}$, ${\bm\Psi _{n1}}$ and $\bm\Lambda_{n1}$ are the top-left $s_n\times s_n$ submatrix of $\bm\Psi_n$ and $\bm\Lambda_{n}$, and ${\bm G_{n1}} = {\bm \Psi_{n1}} - {\bm\mu _{n1}}\bm\mu _{n1}^T$. Let ${\bm\xi _i} = {\bm m_i} - {\bm\mu _{n1}}$ with ${\bm m_i} = {\left( {{m_1}\left( {{X_{i1}}} \right),...,{m_{{s_n}}}\left( {{X_{i{s_n}}}} \right)} \right)^T}$, ${\bm{\tilde \eta }_i} = {\left( {{{\tilde \eta }_{i1}},...,{{\tilde \eta }_{i{s_n}}}} \right)^T}$ with ${{\tilde \eta }_{ij}}  = \sum\limits_{k = 1}^{{s_n}} {{w_{ok}}{\psi _{\tau _k}} \left({\eta _{ik}}\right){\beta _{jk}}\left( {{X_{ik}}} \right)} $ and ${\beta _{jk}}\left( {{x_k}} \right)= E\left[ {{m_j}\left( {{X_{ij}}} \right)\left| {{X_{ik}} = {x_k}} \right.} \right]$. Obviously, the mean of $ { {\left( {{\bm\xi _i}{\psi _\tau }\left( {{\eta _i}} \right) - {f_\eta }\left( 0 \right){\bm{\tilde \eta }_i}} \right)} }$ is zero, and we define $\bm\Upsilon _{n1}=Var\left( {{n^{ - {1 \mathord{\left/
 {\vphantom {1 2}} \right.
 \kern-\nulldelimiterspace} 2}}}\sum\limits_{i = 1}^n {\left( {{\bm\xi _i}{\psi _\tau }\left( {{\eta _i}} \right) - {f_\eta }\left( 0 \right){\bm{\tilde \eta }_i}} \right)} } \right)$. Let ${\bm c_n} = {\left( {{{\dot p}_\lambda }\left( {\left| {{w_{o1}}} \right|} \right){\mathop{\rm sgn}} \left( {{w_{o1}}} \right),...,{{\dot p}_\lambda }\left( {\left| {{w_{o{s_n}}}} \right|} \right){\mathop{\rm sgn}} \left( {{w_{o{s_n}}}} \right)} \right)^T}$, ${\bm\Gamma _n} = diag\left\{ {{{\ddot p}_\lambda }\left( {\left| {{w_{o1}}} \right|} \right),...,{{\ddot p}_\lambda }\left( {\left| {{w_{o{s_n}}}} \right|} \right)} \right\}$ and ${\mathcal{X}_j}=(X_{1j},...,X_{nj})^T$, where ${{\mathop{\rm sgn}} \left(\cdot \right)}$ is the sign function. Define ${\mu _k} = \int {{u^k}} K\left( u \right)du$ and ${\nu _k} = \int {{u^k}} {K^2}\left( u \right)du$. In the following theorems, we give the asymptotic theories of $\hat m_j(\cdot)$ and $ {\bm{{\rm{\hat w}}}_{n}}$.

\noindent \textbf{Theorem 1.} Suppose that $x$ is an interior point of the support of $f_j(\cdot)$. Under the regularity conditions (C1)--(C4), if $h_j\rightarrow 0$ and $nh_j\rightarrow\infty$ for $j=1,...,p_n$, then the asymptotic conditional bias and variance of the local linear estimator $\hat m_j(x)$ are given by
\[\begin{array}{l}
 bias\left( {{{\hat m}_j}\left( x \right)\left| {{\mathcal{X}_j}} \right.} \right) = \frac{1}{2}{{\ddot m}_j}\left( x \right){\mu _2}h_j^2 + {o_p}\left( {h_j^2} \right) ,\\
 Var\left( {{{\hat m}_j}\left( x \right)\left| {{\mathcal{X}_j}} \right.} \right) = \frac{1}{{n{h_j}{f_j}\left( x \right)}}\frac{{\tau_j \left( {1 - \tau_j } \right){\nu _0}}}{{{{\left( {{f_{{\eta _j}}}\left( 0 \right)} \right)}^2}}} + {o_p}\left( {\frac{1}{{n{h_j}}}} \right). \\
 \end{array}\]
Furthermore, conditioning on ${\mathcal{X}_j}$, we have
\[\sqrt {n{h_j}} \left\{ {{{\hat m}_j}\left( x \right) - {m_j}\left( x \right) - \frac{1}{2}{{\ddot m}_j}\left( x \right){\mu _2}h_j^2} \right\}\mathop  \to \limits^d N\left( {0,\frac{1}{{{f_j}\left( x \right)}}\frac{{\tau_j \left( {1 - \tau_j } \right){\nu _0}}}{{{{\left( {{f_{{\eta _j}}}\left( 0 \right)} \right)}^2}}}} \right)\]
for $j=1,...,p_n$, where $\mathop  \to \limits^d$ stands for convergence in distribution.

\noindent \emph{Remark 1.} Theorem 1 shows that the proposed nonparametric
estimate $\hat m_j(\cdot)$ is $\sqrt {n{h_j}}$ consistent and enjoys a asymptotically normal distribution.

\noindent \textbf{Theorem 2.} Under conditions (C1)--(C9), together with $np_n\left( {h_j^2 + {{\left( {n{h_j}} \right)}^{{{ - 1} \mathord{\left/
 {\vphantom {{ - 1} 2}} \right.
 \kern-\nulldelimiterspace} 2}}}} \right) = o\left( 1 \right)$ for $j=1,...,p_n$ and ${{p_n^2} \mathord{\left/
 {\vphantom {{p_n^2} {n \to 0}}} \right.
 \kern-\nulldelimiterspace} {n \to 0}}$, if $\lambda  \to 0,\sqrt {{n \mathord{\left/
 {\vphantom {n {{p_n}}}} \right.
 \kern-\nulldelimiterspace} {{p_n}}}} \lambda  \to \infty $, then $\forall \bm e_n\in \mathcal{R}^{s_n}$ such that $\left\|\bm e_n\right\|=1$, we have

(i) there exists a local minimizer ${\bm{{\rm{\hat w}}}_{n}}$ of the objective function $\mathcal{Q}_n\left( {\bm{{\rm{w}}}_{n}}\right)$ defined in (2) such that
$\left\| { {\bm{{\rm{\hat w}}}_{n}} -  {\bm{{\rm{ w}}}_{o}}} \right\| = {O_p}\left( {\sqrt {{p_n}} \left( {{n^{{{ - 1} \mathord{\left/
 {\vphantom {{ - 1} 2}} \right.
 \kern-\nulldelimiterspace} 2}}} + {a_n}} \right)} \right)$;

 (ii) ${\bm{{\rm{\hat w}}}_{n}(2)}=\bm 0$ with probability approaching one;

 (iii) $\sqrt n \bm e_n^T\bm \Upsilon_{n1}^{{{ - 1} \mathord{\left/
 {\vphantom {{ - 1} 2}} \right.
 \kern-\nulldelimiterspace} 2}}\left( {{f_\eta }\left( 0 \right)\bm G_{n1} + {\bm\Gamma _n}} \right)\left[ {{\bm{{\rm{\hat w}}}_{n}(1)} - {\bm{{\rm{ w}}}_{o}(1)} + {{\left( {{f_\eta }\left( 0 \right)\bm G_{n1} + {\bm\Gamma _n}} \right)}^{ - 1}}{\bm c_n}} \right]\mathop  \to \limits^d N\left( 0,1\right)$.

\noindent \emph{Remark 2.} Theorem 2 indicates that the estimate of the optimal weight $\bm{{\rm{ w}}}$ is still consistent although the dimension of predictor goes to infinite. Meanwhile, it also shows that the proposed estimate ${\bm{{\rm{\hat w}}}_{n}}$  enjoys well-known properties in high dimensional variable selection such as the sparsity and oracle property.

\section{Numerical studies}
We investigate the performance of the proposed approach by three simulation examples and an empirical application. In our numerical studies, we set the kernel function $K(\cdot)$ as the Epanechnikov kernel, namely, $K\left( u \right) = 0.75{\left( {1 - {u^2}} \right)_ + }$. Bandwidth selection is crucial in local smoothing since it governs the curvature of the fitted function. Similar to Kai et al. (2011), we use the following formula to choose the bandwidth $h=h_{ls}\times\left\{ {\tau(1-\tau)/f(\Phi^{-1}(\tau))} \right\} ^{1/5}$, where $h_{ls}$ is the selected optimal bandwidth for least squares, and $f(\cdot)$ and $\Phi(\cdot)$ represent the density function and distribution function of standard normal distribution, respectively. The rule of thumb is used to select the bandwidth $h_{ls}$. In addition, the tuning parameter $\lambda$ in the proposed penalized procedure plays an important role. Lian (2012) had proved that the Schwarz
information criterion (SIC) is a consistent variable selection criterion under the framework of fixed dimension. In this paper, we select $\lambda$ by minimizing the following modified SIC criterion (MSIC)
\begin{equation}
{\rm {MSIC}}\left( {{\lambda}} \right) = \log \left(\mathcal{Q}_n\left( {\bm{{\hat{\rm{w}}}}_{n}}\right) \right) + {{df C_n\log \left( n \right)} \mathord{\left/
 {\vphantom {{df{C_n}\log \left( n \right)} {\left( {2n} \right)}}} \right.
 \kern-\nulldelimiterspace} {\left( {2n} \right)}},
 \end{equation}
where ${\bm{{\hat{\rm{w}}}}_{n}} $ is the estimated model weight vector for a given $\lambda$, $df$ is the number of nonzero coefficients in ${\bm{{\hat{\rm{w}}}}_{n}} $ and $C_n$ diverges with $n$. For example, the MSIC criterion reduces to tradition SIC criterion Lian (2012) when $C_n=1$, and the MSIC criterion is more suitable for high dimensional data if $C_n$ is selected as $\log \left( p_n \right)$.

In order to investigate the superiority of the proposed method, we consider the following methods:
(1) the proposed semiparametric model average quantile prediction (without SCAD penalty, denoted as SMAQP), (2)the proposed penalized semiparametric model average quantile prediction (with SCAD penalty, denoted as PSMAQP), (3) semiparametric model average mean prediction proposed by Li et al. (2015) (without SCAD penalty, denoted as SMAMP), (4) penalized semiparametric model average mean prediction proposed by Chen et al. (2018) (with SCAD penalty, denoted as PSMAMP). SMAMP and PSMAMP aim to forecast the conditional mean function $E\left( Y|\bm X \right)$, and detailed descriptions about the two methods can refer to the section 3 of Li et al. (2015) and subsection 2.1 of Chen et al. (2018). The tuning parameter $\lambda$ involved in PSMAMP is chosen by the cross-validation according to the advice of Chen et al. (2018), and the R package ``ncvreg'' can be used to obtain the penalized estimator PSMAMP.

 \subsection{Simulation studies}
 In all simulation examples, the sample size $n$ consists of a training set of size $n_{tr}$ and a testing set of size $n_{te}$, namely, $n=n_{tr}+n_{te}$.

\textbf{Example 1}. For a clear comparison, we adopt similar settings used in Chen et al. (2018) and generate the random samples from the following model
\begin{equation}
{Y_i} = {{\rm{m}}_{\rm{1}}}\left( {{X_{i1}}} \right) + {{\rm{m}}_2}\left( {{X_{i2}}} \right) + {{\rm{m}}_3}\left( {{X_{i3}}} \right) + {{\rm{m}}_4}\left( {{X_{i4}}} \right) + {\varepsilon _i}, i=1,...,n,
 \end{equation}
where ${{\rm{m}}_{\rm{1}}}\left( u \right)={\rm{-sin}}(2u)$, ${{\rm{m}}_{\rm{2}}}\left( u \right)=u^2-25/12$, ${{\rm{m}}_{\rm{3}}}\left( u \right)=u$ and ${{\rm{m}}_{\rm{4}}}\left( u \right)=e^{-u}-\frac{2}{5}{\rm{sinh}}(5/2)$. We fix $n_{te}=100$ and consider $n_{tr}=200$ and $400$ for example 1. The covariates $\bm X_i=\left(X_{i1},...,X_{ip_n}\right)^T$ are independently drawn from $U\left(-2.5,2.5\right)$, and we set the dimension of covariates as $p_n=[n_{tr}^{1/2}]$ which satisfies the theoretical condition ${{p_n^2} \mathord{\left/
 {\vphantom {{p_n^2} {n \to 0}}} \right. \kern-\nulldelimiterspace} {n \to 0}}$, where $[s]$ stands for the largest integer not greater than $s$. Obviously, the first four variables make a significant contribution to estimating the joint multivariate quantile function $m_{\tau}\left(\bm X\right)$, while the rest are not. Therefore, we have reasons to believe that the first four model weights are nonzero and the rest are zero. Please note that the model average component $m_j $ given in section 2 is different from ${{\rm{m}}_j}$ reported in model (4) for $j=1,...,p_n$. Our mission is to achieve the goal of accurately predicting the conditional quantile function $m\left( \cdot\right)$, so we are not attempting to estimate ${{\rm{m}}_j}$ in this paper.

In order to examine the robustness of the proposed procedure, we consider the following three different error distributions of $\varepsilon_i$: standard normal distribution (SN), $t$-distribution with 3 degrees of freedom ($t_3$), contaminated normal distribution ($\rm {MN}(\rho, \sigma_1,\sigma_2)$) representing a mixture of $\rm {N}(0,1)$ and $\rm {N}(0, 10^2)$ with weights $0.95$ and $0.05$ respectively. In addition, four criteria are adopted to evaluate the performance of proposed approach. Firstly, ``C'', ``IC'' and ``CF'' are considered to examine variable selection performance, where ``C'' represents the average number of zero coefficients in the model weight vector that are correctly estimated to be zero; ``IC'' represents the average number of nonzero coefficients in the model weight vector that are incorrectly estimated to be zero and ``CF'' represents the proportion of correctly fitted models (``correctly fit'' means that the estimation procedure correctly chooses all significant
components from the model weight vector). Secondly, the mean prediction error (MPE) is used to measure accuracy of prediction, which is defined as $\sum\nolimits_{i \in \mathcal{I}} {{\rho _\tau }\left( {{Y_i} - {{\hat Y}_i}} \right)} /\left| \mathcal{I} \right|$, where $\mathcal{I}$ stands for an index set of either the training sample or the testing sample.

\tabcolsep=11pt
\begin{table}\scriptsize
\caption{Simulation results of C, IC, CF, MPE and their standard deviations (in parenthesis) for $\tau=0.5$ in example 1.}
\begin{tabular}{ccccccccc} \noalign{\smallskip}\hline
\multicolumn{1}{c}{$n_{tr}$}
&\multicolumn{1}{c}{Error}
&\multicolumn{1}{c}{method}
&\multicolumn{1}{c}{C}
&\multicolumn{1}{c}{IC}
&\multicolumn{1}{c}{CF}
&\multicolumn{1}{c}{In-Sample Error}
&\multicolumn{1}{c}{Out-of-Sample Error}
\\
\hline
200&N(0,1)&SMAMP&--&--&--& 0.459 (0.031)& 0.577 (0.058)\\
&&PSMAMP& 7.362 &0.018& 0.258&0.467 (0.032)& 0.566 (0.058)\\
&&SMAQP& --&--&--&0.462 (0.036)& 0.580 (0.062)\\
&&PSMAQP&9.930 &0.108 &0.836&0.487 (0.038)&0.557 (0.061) \\
&$t_3$&SMAMP&--&--&--& 0.601 (0.065) &0.748 (0.085)\\
&&PSMAMP&6.118& 0.030& 0.096 &0.605 (0.063)& 0.730 (0.085)\\
&&SMAQP&--&--&--&0.592 (0.047) &0.734 (0.081) \\
&&PSMAQP&9.956& 0.206 &0.760&0.625 (0.051)& 0.705 (0.077) \\
&MN&SMAMP&--&--&--& 0.702 (0.117) &0.869 (0.141)\\
&&PSMAMP&4.372& 0.048& 0.030& 0.697 (0.113) &0.846 (0.145)\\
&&SMAQP&--&--&--& 0.637 (0.084) &0.765 (0.122)\\
&&PSMAQP&9.966 &0.222& 0.754& 0.668 (0.087) &0.738 (0.121)\\
\hline
400&N(0,1)&SMAMP&--&--&--&0.445 (0.021)& 0.524 (0.046) \\
&&PSMAMP&13.63 &0.006 &0.414& 0.453 (0.022) &0.513 (0.045)\\
&&SMAQP&--&--&--&0.439 (0.023)& 0.512 (0.045) \\
&&PSMAQP&15.98& 0.006 &0.976&0.457 (0.024)& 0.494 (0.043) \\
&$t_3$&SMAMP&--&--&--&0.597 (0.052)& 0.694 (0.081) \\
&&PSMAMP&11.38& 0.002& 0.162&0.599 (0.042)& 0.673 (0.078) \\
&&SMAQP&--&--&--&0.580 (0.036) &0.664 (0.072) \\
&&PSMAQP&15.99& 0.012 &0.978& 0.603 (0.037) &0.640 (0.072)\\
&MN&SMAMP&--&--&--&0.676 (0.078) &0.789 (0.116) \\
&&PSMAMP&8.346 &0.006& 0.040&0.667 (0.075) &0.761 (0.120) \\
&&SMAQP&--&--&--&0.616 (0.058)& 0.693 (0.106) \\
&&PSMAQP& 15.99& 0.012& 0.984&0.636 (0.058)& 0.674 (0.105)\\
\hline
\end{tabular}
\\
Notation: To make this a fair comparison, we consider $\tau=0.5$ for SMAQP and PSMAQP in this table. In addition, the number of zero components of model weight vector is 10 for $n_{tr}=200$ and 16 for $n_{tr}=400$.
\end{table}

\textbf{Example 2}. In this example, similar to Huang and Li (2018), we generate the random samples from the following model
\begin{equation}
{Y_i} = 3{{\rm{m}}_{\rm{1}}}\left( {{X_{i1}}} \right) + 3{{\rm{m}}_2}\left( {{X_{i2}}} \right) + 2{{\rm{m}}_3}\left( {{X_{i3}}} \right) + 2{{\rm{m}}_4}\left( {{X_{i4}}} \right) + \sqrt{1.74}{\varepsilon _i}, i=1,...,n,
 \end{equation}
where ${{\rm{m}}_{\rm{1}}}\left( u \right)=(2u)$, ${{\rm{m}}_{\rm{2}}}\left( u \right)=(2u-1 )^2 $, ${{\rm{m}}_{\rm{3}}}\left( u \right)={\rm{sin}}(2\pi u)/(2-{\rm{sin}}(2\pi u))$ and ${{\rm{m}}_{\rm{4}}}\left( u \right)=0.1{\rm{sin}}(2\pi u)+0.2{\rm{cos}}(2\pi u)+0.3{\rm{sin}}(2\pi u)^2+0.4{\rm{cos}}(2\pi u)^3+0.5{\rm{sin}}(2\pi u)^3$. The covariates $\bm X_i=\left(X_{i1},...,X_{ip_n}\right)^T$ are simulated by ${X_{ij}} = \left( {{W_{ij}} + t{U_i}} \right)/\left( {1 + t} \right)$ for $t=1$ and $j=1,...,p_n$, where $W_{ij}$ and ${U_i}$ are independently drawn from $U(0, 1)$ and $p_n=[n_{tr}^{1/2}]$. We also fix $n_{te}=100$ and consider $n_{tr}=400$ and $800$ for example 2. Other settings are the same as that in example 1.

It is easy to find that the conditional mean function $E\left( Y|\bm X \right)$ is equal to the conditional quantile function $m_{\tau}\left(\bm X\right)$ for $\tau=0.5$. Thus, we can compare mean prediction approaches (SMAMP and PSMAMP) with quantile prediction approaches (SMAQP and PSMAQP) at $\tau=0.5$. The MPE criterion is reduced to $\sum\nolimits_{i \in \mathcal{I}} { \frac{1}{2}\left| {{{Y_i} - {{\hat Y}_i}}} \right|} /\left| \mathcal{I} \right|$ for $\tau=0.5$, and thus this criterion also can be used to assess the prediction performance of mean prediction approaches. The corresponding results of mean prediction approaches (SMAMP and PSMAMP) and quantile prediction approaches (SMAQP and PSMAQP) at $\tau=0.5$ are reported in Tables 1 and 3. We can obtain the following findings. Firstly, the values in the column labeled ``C'' gradually tend to the true number of zero components with the training sample size increasing. The CF values are very close to one for a large training sample size (e.g. $n_{tr}=400$), which shows that the proposed penalized procedure can consistently select significant components in weight vector. However, the existing mean prediction approach PSMAMP performs badly due to lower CF values. Secondly, unpenalized methods always has smaller in-sample MPE than the penalized methods's, but it does not hold true for out-of-sample MPE. For heavy-tailed distributions $t_3$ and contaminated distribution MN, it is not hard to find that our proposed penalized method PSMAQP is best in terms of prediction accuracy among all methods. Meanwhile, there is little difference for PSMAMP and PSMAQP under the normal error distribution. Thirdly, Tables 2 and 4 give the simulation results of SMAQP and PSMAQP at $\tau=0.75$. The results also show that PSMAQP has better prediction performance.

\tabcolsep=11pt
\begin{table}\scriptsize
\caption{Simulation results of C, IC, CF, MPE and their standard deviations (in parenthesis) for $\tau=0.75$ in example 1.}
\begin{tabular}{ccccccccc} \noalign{\smallskip}\hline
\multicolumn{1}{c}{$n_{tr}$}
&\multicolumn{1}{c}{Error}
&\multicolumn{1}{c}{method}
&\multicolumn{1}{c}{C}
&\multicolumn{1}{c}{IC}
&\multicolumn{1}{c}{CF}
&\multicolumn{1}{c}{In-Sample Error}
&\multicolumn{1}{c}{Out-of-Sample Error}
\\
\hline
200&N(0,1)&SMAQP&--&--&--& 0.917 (0.099)& 1.032 (0.163)\\
&&PSMAQP&9.896 &0.562& 0.468& 0.961 (0.108) &1.027 (0.164)\\
&$t_3$&SMAQP& --&--&--&1.087 (0.110) &1.239 (0.170)\\
&&PSMAQP& 9.866& 0.768 &0.344&1.139 (0.130) &1.229 (0.177)\\
&MN&SMAQP&--&--&--& 1.108 (0.128) &1.256 (0.213)\\
&&PSMAQP& 9.918& 0.848 &0.338&1.168 (0.148)&1.253 (0.215)\\
\hline
400&N(0,1)&SMAQP&--&--&--& 0.850 (0.075)& 0.921 (0.119)\\
&&PSMAQP& 15.97  &0.134&  0.844&0.875 (0.077)& 0.909 (0.121)\\
&$t_3$&SMAQP&--&--&--& 1.035 (0.085) &1.129 (0.148)\\
&&PSMAQP& 15.95 & 0.224 & 0.740&1.066 (0.088)& 1.110 (0.147)\\
&MN&SMAQP& --&--&--&1.054 (0.096) &1.147 (0.181)\\
&&PSMAQP& 15.97&  0.208&  0.776&1.081 (0.100)& 1.127 (0.178)\\
\hline
\end{tabular}
\\
Notation: The number of zero components of model weight vector is 10 for $n_{tr}=200$ and 16 for $n_{tr}=400$.
\end{table}

\tabcolsep=11pt
\begin{table}\scriptsize
\caption{Simulation results of C, IC, CF, MPE and their standard deviations (in parenthesis) for $\tau=0.5$ and $t=1$ in example 2.}
\begin{tabular}{ccccccccc} \noalign{\smallskip}\hline
\multicolumn{1}{c}{$n_{tr}$}
&\multicolumn{1}{c}{Error}
&\multicolumn{1}{c}{method}
&\multicolumn{1}{c}{C}
&\multicolumn{1}{c}{IC}
&\multicolumn{1}{c}{CF}
&\multicolumn{1}{c}{In-Sample Error}
&\multicolumn{1}{c}{Out-of-Sample Error}
\\
\hline
400&N(0,1)&SMAMP&--&--&--&0.507 (0.021) &0.566 (0.045)\\
&&PSMAMP&11.75 &0.000 &0.212& 0.515 (0.022) &0.555 (0.045)\\
&&SMAQP& --&--&--&0.501 (0.020) &0.594 (0.047)\\
&&PSMAQP&15.97 &0.006& 0.966&0.530 (0.022) &0.566 (0.045)\\
&$t_3$&SMAMP&--&--&--&0.724 (0.054)& 0.810 (0.105)\\
&&PSMAMP& 7.750& 0.006& 0.026&0.725 (0.052) &0.790 (0.107)\\
&&SMAQP&--&--&--& 0.692 (0.042) &0.802 (0.104)\\
&&PSMAQP&15.98 &0.048& 0.936&0.728 (0.044)& 0.769 (0.103)\\
&MN&SMAMP&--&--&--&0.837 (0.105) &0.932 (0.154)\\
&&PSMAMP& 5.740 &0.012& 0.000&0.826 (0.104)& 0.908 (0.156)\\
&&SMAQP& --&--&--&0.740 (0.072) &0.839 (0.151)\\
&&PSMAQP&15.99 &0.044& 0.946&0.772 (0.073)& 0.808 (0.151)\\
\hline
800&N(0,1)&SMAMP& --&--&--&0.511 (0.014)& 0.553 (0.044)\\
&&PSMAMP& 18.61  &  0& 0.276&0.519 (0.015) &0.544 (0.043)\\
&&SMAQP& --&--&--& 0.508 (0.014)& 0.567 (0.045)\\
&&PSMAQP&23.99& 0 &0.988&0.528 (0.014) &0.546 (0.043)\\
&$t_3$&SMAMP&--&--&--&0.728 (0.038)& 0.795 (0.089)\\
&&PSMAMP& 13.89&    0& 0.036&0.727 (0.035)& 0.774 (0.090)\\
&&SMAQP&--&--&--& 0.704 (0.030)& 0.784 (0.088)\\
&&PSMAQP&23.99&   0& 0.998&0.730 (0.031)& 0.758 (0.087)\\
&MN&SMAMP&--&--&--&0.817 (0.066)& 0.884 (0.140)\\
&&PSMAMP&9.662&   0& 0.008& 0.803 (0.067)& 0.856 (0.146)\\
&&SMAQP&--&--&--&0.744 (0.048)& 0.810 (0.137) \\
&&PSMAQP&24.00&   0& 1.000&0.767 (0.048)& 0.787 (0.136)\\
\hline
\end{tabular}
\\
Notation: To make this a fair comparison, we consider $\tau=0.5$ for RSMAP and PRSMAP in this table. In addition, the number of zero components of model weight vector is 16 for $n_{tr}=400$ and 24 for $n_{tr}=800$.
\end{table}

\tabcolsep=11pt
\begin{table}\scriptsize
\caption{Simulation results of C, IC, CF, MPE and their standard deviations (in parenthesis) for $\tau=0.75$ and $t=1$ in example 2.}
\begin{tabular}{ccccccccc} \noalign{\smallskip}\hline
\multicolumn{1}{c}{$n_{tr}$}
&\multicolumn{1}{c}{Error}
&\multicolumn{1}{c}{method}
&\multicolumn{1}{c}{C}
&\multicolumn{1}{c}{IC}
&\multicolumn{1}{c}{CF}
&\multicolumn{1}{c}{In-Sample Error}
&\multicolumn{1}{c}{Out-of-Sample Error}
\\
\hline
400&N(0,1)&SMAQP& --&--&--&0.848 (0.042)& 0.938 (0.095)\\
&&PSMAQP&15.92 & 0.010 & 0.912&0.882 (0.043)& 0.922 (0.094)\\
&$t_3$&SMAQP& --&--&--&1.111 (0.070)& 1.220 (0.139)\\
&&PSMAQP&15.91  &0.110&  0.828&1.147 (0.073) &1.188 (0.133)\\
&MN&SMAQP& --&--&--&1.121 (0.092)& 1.204 (0.195)\\
&&PSMAQP&15.93 & 0.094 & 0.852&1.160 (0.094)& 1.184 (0.193)\\

800&N(0,1)&SMAQP&--&--&--&0.853 (0.029)& 0.912 (0.083)\\
&&PSMAQP&23.97 & 0.000&  0.968&0.879 (0.030)& 0.901 (0.081)\\
&$t_3$&SMAQP&--&--&--&1.114 (0.049)& 1.179 (0.150)\\
&&PSMAQP&23.97 & 0.002 & 0.968&1.136 (0.048) &1.151 (0.143)\\
&MN&SMAQP&--&--&--&1.116 (0.065)& 1.165 (0.179)\\
&&PSMAQP&23.98  &0.000&  0.984& 1.142 (0.065)& 1.146 (0.174)\\
\hline
\end{tabular}
\\
Notation: The number of zero components of model weight vector is 16 for $n_{tr}=400$ and 24 for $n_{tr}=800$.
\end{table}

\tabcolsep=12pt
\begin{table}\scriptsize
\caption{Simulation results of C, IC, CF, MEE and their standard deviations (in parenthesis) in example 3.}
\begin{tabular}{ccccccccc} \noalign{\smallskip}\hline
\multicolumn{1}{c}{$n_{tr}$}
&\multicolumn{1}{c}{$\tau$}
&\multicolumn{1}{c}{method}
&\multicolumn{1}{c}{C}
&\multicolumn{1}{c}{IC}
&\multicolumn{1}{c}{CF}
&\multicolumn{1}{c}{In-Sample Error}
&\multicolumn{1}{c}{Out-of-Sample Error}
\\
\hline
400&0.5&SMAQP&--&--&--& 0.328 (0.038)& 0.341 (0.047)\\
&&PSMAQP& 14.96&  0.350&  0.652&0.258 (0.049)& 0.266 (0.056)\\
&0.75&SMAQP&--&--&--&0.383 (0.047)& 0.401 (0.060)\\
&&PSMAQP&14.90&  0.236 & 0.702&0.301 (0.058)& 0.310 (0.066)\\

800&0.5&SMAQP&--&--&--&0.266 (0.027)& 0.273 (0.036)\\
&&PSMAQP&22.99&  0.014 & 0.974&0.177 (0.029)& 0.180 (0.033)\\
&0.75&SMAQP&--&--&--&0.316 (0.033)& 0.325 (0.046)\\
&&PSMAQP&22.95&  0.022 & 0.932&0.217 (0.037) &0.221 (0.042)\\
\hline
\end{tabular}
\\
Notation: The number of zero components of model weight vector is 15 for $n_{tr}=400$ and 23 for $n_{tr}=800$.
\end{table}

\textbf{Example 3}. The conditional quantile function $ Q_{\tau}\left( {Y|\bm X} \right)$ is considered as
\begin{equation}
\begin{array}{l}
m_{\tau}\left(\bm X\right)=Q_{\tau}\left( {Y|\bm X} \right)= 1+{\Phi ^{ - 1}}\left( \tau  \right) + 2{X_{i1}} + 3X_{i2}^2 - \log \left( {1 - {X_{i3}}} \right) + {\Phi ^{ - 1}}\left( {{X_{i4}}} \right) +  \\
~~~~~~~~~~~~~~~~~~~~~~~~~~~~~~ \left( {1 + {\Phi ^{ - 1}}\left( \tau  \right) - \log \left( {1 - \tau } \right)} \right){X_{i5}}, i=1,...,n, \\
 \end{array}
\end{equation}
where $\Phi$ is the standard normal distribution function, $\bm X_i=\left(X_{i1},...,X_{ip_n}\right)^T$ are independently drawn from $U\left(0,1\right)$ and $1+{\Phi ^{ - 1}}\left( \tau  \right)$ can be regarded as the intercept. We fix $n_{te}=100$ and consider $n_{tr}=400$ and $800$ for example 3 and $p_n=[n_{tr}^{1/2}]$. Obviously, the fifth covariate's coefficient varies with $\tau$, and only the first five predictors are significant for predicting $Q_{\tau}\left( {Y|\bm X} \right)$. The first two examples come from the nonparametric additive model, but the proposed approach do not need any model assumption. Thus we aim to confirm that our method is model free in this example. To assess the estimation accuracy of $\hat m_{\tau}\left(\bm X\right)$, we consider the mean estimation error (MEE) defined as ${\rm {MEE}}=\sum\nolimits_{i \in \mathcal{I}} { \frac{1}{2}\left| {{ m_{\tau}\left(\bm X\right) - \hat m_{\tau}\left(\bm X\right)}} \right|} /\left| \mathcal{I} \right|$ in this example. Table 5 lists the simulation results which show that the proposed PSMAQP performs well for different quantiles.

Overall, the proposed model free procedure PSMAQP is competitive when compared with the existing methods, and its finite sample performances are satisfactory.

\subsection{An application}
In this section, we apply our proposed method to analyze the body fat dataset (Johnson, 1996), which is available from \emph{http://lib.stat.cmu.edu/datasets/bodyfat}. This dataset consists of 252 observations without missing values. The purpose of studying this dataset is to predict the percentage of body fat according to various body circumference measurements. Thus, the percentage of body fat is taken as the response variable and other body circumference measurements are regarded as the predictors. Brief descriptions and marginal Pearson correlations of 14 variables are summarized in Table 6. More details can refer to Johnson (1996). Before employing prediction methods, we take the logarithm transformation for all predictors.

To evaluate the predictive performance of various methods, the data is split into two parts. One part including $n_{tr}$ observations is used as a training data set to estimate the weight vector ${\bm{{\rm{w}}}_{n}}$ and the marginal quantile functions $m_j(\cdot), j=1,...,p_n$, and the other part including $n_{te}$ observations is considered as a testing data set to evaluate the predictive ability of various methods. In this real data analysis, we consider $n_{tr}=150$ and $200$, and $n_{te}=n-n_{tr}=252-n_{tr}$.

Table 7 reports the in-sample and out-of-sample mean prediction errors (MPE) and the corresponding sample standard deviations (SD) over 500 random partitions. Firstly, for $\tau=0.5$ and in-sample performance, it is easy to see that SMAQP performs best among four approaches. For out-of-sample performance, one can see clearly that the proposed penalized approach PSMAQP has smallest MPE and SD for different settings, which shows that our proposed method has better predictive ability. Secondly, for $\tau=0.25$ and $0.75$, we can see that PSMAQP always performs better than SMAQP in terms of out-of-sample performance.

To investigate the estimated weights, we list the estimated weights at $\tau=0.5$ and their standard deviations (in brackets) calculated by the bootstrap resampling method (Horowitz, 1998). Obviously, the weights for the penalized prediction methods (PSMAMP and PSMAQP) are relatively sparse with much smaller standard deviations than the unpenalized prediction methods (SMAMP and SMAQP). Meanwhile, it is not hard to find that PSMAQP is most efficient among all methods due to the smallest standard deviations. In addition, for PSMAMP, only the sixth predictor ($X_6$) is chose as the significant variable whose marginal quantile function has significant influence on estimating $m_{\tau}(\bm X)$. However, PSMAQP selects five predictors (including $X_1$, $X_4$, $X_6$, $X_9$ and $X_{13}$) as the significant variables. In summary, our proposed model averaging procedure generally works
well and outperforms other existing methods.

\tabcolsep=12pt
\begin{table}\scriptsize
\caption{\leftline{Regressors for the body fat dataset.}}
\begin{tabular}{ccccccccc} \noalign{\smallskip}\hline
\multicolumn{1}{c}{Variable}
&\multicolumn{1}{c}{Name}
&\multicolumn{1}{c}{Description}
&\multicolumn{1}{c}{Correlation with $Y$}
\\
\hline
$X_1$&Age&Age (years)&0.2921\\
$X_2$&Weight& Weight(lbs)&0.6287\\
$X_3$&Height&Height (inches)&-0.0990\\
$X_4$&Neck& Neck circumference (cm)&0.4905\\
$X_5$&Chest& Chest circumference (cm)&0.7051\\
$X_6$&Abdomen&Abdomen 2 circumference (cm)&0.8218\\
$X_7$&Hip&Hip circumference (cm)&0.6365\\
$X_8$&Thigh &Thigh circumference (cm)&0.5680\\
$X_9$&Knee &Knee circumference (cm)&0.5102\\
$X_{10}$&Ankle &Ankle circumference (cm)&0.2796\\
$X_{11}$&Biceps &Biceps (extended) circumference (cm)&0.5000\\
$X_{12}$&Forearm &Forearm circumference (cm)&0.3584\\
$X_{13}$&Wrist &Wrist circumference (cm)& 0.3447\\
$Y$&Body Fat($\%$)&Percent body fat using Siri's equation: 495/Density - 450&1.0000\\
\hline
\end{tabular}
\end{table}

\tabcolsep=9pt
\begin{table}\scriptsize
\caption{\leftline{Prediction results ($\times 10^{-2}$) for analysis of the body fat dataset.}}
\label{table1}
\begin{tabular}{cccccccccccccc} \noalign{\smallskip}\hline
&\multirow{3}{*}{method}
&\multicolumn{5}{c}{$n_{tr}=150$}
&\multicolumn{5}{c}{$n_{tr}=200$}
\\
\cline{3-7}
\cline{9-13}
&&\multicolumn{2}{c}{In-Sample}
&&\multicolumn{2}{c}{Out-of-Sample}
&&\multicolumn{2}{c}{In-Sample}
&&\multicolumn{2}{c}{Out-of-Sample}
\\
\cline{3-4}
\cline{6-7}
\cline{9-10}
\cline{12-13}
&&\multicolumn{1}{c}{MPE}&\multicolumn{1}{c}{SD}
&&\multicolumn{1}{c}{MPE}&\multicolumn{1}{c}{SD}
&&\multicolumn{1}{c}{MPE}&\multicolumn{1}{c}{SD}
&&\multicolumn{1}{c}{MPE}&\multicolumn{1}{c}{SD}
\\
\hline
$\tau=0.5$&SMAMP&1.622 &0.069 &&2.116& 0.691&&1.653 &0.043&& 2.242& 1.402\\
&PSMAMP&1.763 &0.136&& 2.117 &0.544&&1.829 &0.105&& 2.314& 1.357\\
&SMAQP&1.593 &0.067&& 1.940& 0.143&&1.634& 0.043&& 1.890& 0.174\\
&PSMAQP&1.683 &0.081&& 1.896 &0.125&&1.696& 0.055&& 1.858& 0.166\\
\hline
$\tau=0.25$&SMAQP&2.666& 0.151 &&2.955& 0.338&&2.687 &0.098&& 2.857 &0.363\\
&PSMAQP&2.726 &0.165&& 2.914 &0.314&&2.742& 0.089&& 2.853& 0.348\\
\hline
$\tau=0.75$&SMAQP&2.768& 0.169&& 3.118 &0.383&&2.838& 0.115 &&3.079 &0.451\\
&PSMAQP&2.811 &0.171&& 3.080& 0.340&&2.859 &0.122 &&3.056 &0.401\\
\hline
\end{tabular}
\end{table}

\tabcolsep=22pt
\begin{table}\scriptsize
\caption{\leftline{Estimated weights and their standard deviations (in brackets) for the body fat study at $\tau=0.5$.}}
\label{table1}
\begin{tabular}{cccccccccccccc} \noalign{\smallskip}\hline
\multicolumn{1}{c}{weight}
&\multicolumn{1}{c}{SMAMP}
&\multicolumn{1}{c}{PSMAMP}
&\multicolumn{1}{c}{SMAQP}
&\multicolumn{1}{c}{PSMAQP}
\\
\hline
$w_0$&-0.056 (12.091)& -0.044 (15.681)& -0.038 (0.584)&  0.063 (0.070)\\
$w_1$&0.660  (0.306)&  0.000  (0.391) & 0.419 (0.189)&  0.314 (0.210)\\
$w_2$&-0.530  (0.546)&  0.000  (0.446)& -0.115 (0.546)&  0.000 (0.231)\\
$w_3$&1.114 (3.454) & 0.000 (2.229)&  0.146 (2.976)&  0.000 (0.192)\\
$w_4$&-0.418  (0.338)&  0.000  (0.397)& -0.330 (0.216)& -0.238 (0.225)\\
$w_5$&-0.055  (0.210)&  0.000  (0.192)&  0.005 (0.248)&  0.000 (0.128)\\
$w_6$&1.643  (0.151)&  1.234  (0.185)&  1.174 (0.185)&  1.187 (0.170)\\
$w_7$&-0.188  (0.291)&  0.000  (0.259)& -0.164 (0.282)&  0.000 (0.190)\\
$w_8$&0.515  (0.251) & 0.000  (0.291)&  0.282 (0.178)&  0.000 (0.188)\\
$w_9$&-0.185  (0.336) & 0.000  (0.264)& -0.195 (0.195)& -0.049 (0.167)\\
$w_{10}$&0.330  (0.752)&  0.000  (0.590)&  0.204 (0.312)&  0.000 (0.156)\\
$w_{11}$&0.292  (0.209)&  0.000  (0.214)&  0.124 (0.166)&  0.000 (0.167)\\
$w_{12}$&0.204  (0.298)&  0.000  (0.245)&  0.343 (0.205)&  0.000 (0.158)\\
$w_{13}$&-2.095  (1.273)&  0.000  (1.519)& -0.696 (0.303)& -0.561 (0.349)\\
\hline
\end{tabular}
\end{table}

\section{Conclusion}

In this paper, we provide a new semiparametric model averaging estimation for forecasting the conditional quantile function $m_{\tau}\left(\bm X\right)$ under the high-dimensional settings. Based on local linear regression, we firstly estimate the individual marginal regression functions $m_j(\cdot)$ by minimizing the local weighted quantile loss function. Then, a penalized quantile regression is developed to select the regressors whose marginal regression functions make significant contribution in estimating the quantile function $m_{\tau}\left(\bm X\right)$. Simulations and empirical example in Section 4 show that the proposed method performs reasonably well in finite samples.

Recently, under the framework of ultra-high dimension setting, Ando and Li (2014) developed a new model averaging approach based on delete-one cross-validation criterion and proved that the proposed method could achieve the lowest possible prediction loss asymptotically.
But they only considered high dimensional parametric model averaging, which may increase the risk of model misspecification. Thus, it is interesting to study semiparametric model averaging estimation for ultra-high dimensional data. Research in these aspects is ongoing.

\noindent \textbf{Acknowledgments}
\\
 This work is supported by the National Social Science Fund of China (Grant No. 17CTJ015).
\\

\noindent \textbf{Appendix}

Let $C$ denote a positive constant that may be different at different place throughout this paper.
Let $u = \sqrt {n{h_j}} \left( {a - {m_j}\left( x \right)} \right),v= {h_j}\sqrt {n{h_j}} \left( {b - {{\dot m}_j}\left( x \right)} \right),{x_{ij}} = {{\left( {{X_{ij}} - x} \right)} \mathord{\left/
 {\vphantom {{\left( {{X_{ij}} - x} \right)} {{h_j}}}} \right.
 \kern-\nulldelimiterspace} {{h_j}}}$, ${K_{ij}} = K\left( {{x_{ij}}} \right)$, ${d_{ij}} = {m_j}\left( {{X_{ij}}} \right) - {m_j}\left( x \right) - {{\dot m}_j}\left( x \right)\left( {{X_{ij}} - x} \right)$ and ${\Delta _{ij}} = {{\left( {u + v{x_{ij}}} \right)} \mathord{\left/
 {\vphantom {{\left( {u + v{x_{ij}}} \right)} {\sqrt {n{h_j}} }}} \right.
 \kern-\nulldelimiterspace} {\sqrt {n{h_j}} }}$. Define $S_{ij}^* = I\left( {{\eta _{ij}} + {d_{ij}} \leq 0} \right) - \tau_j ,\bm W_n^* = {\left( {w_1^*,w_2^*} \right)^T},w_1^* = \frac{1}{{\sqrt {n{h_j}} }}\sum\limits_{i = 1}^n {{K_{ij}}S_{ij}^*} ,w_2^* = \frac{1}{{\sqrt {n{h_j}} }}\sum\limits_{i = 1}^n {{K_{ij}}S_{ij}^*{x_{ij}}} $.

\noindent \textbf{Lemma 1.} Denote $\bm {\hat \theta}  = {\left( {{\hat u},\hat v} \right)^T}$ as the minimizer of (1). Then, under the regularity conditions (C1)--(C4), we have
\[\bm{\hat \theta}  + \frac{1}{{{f_j}\left( x \right)}}{\bm S^{ - 1}}E\left[ {\bm W_n^*\left| {{\mathcal{X}_j}} \right.} \right]\mathop  \to \limits^d N\left( {\bm 0,\frac{1}{{{f_j}\left( x \right)}}{\bm S^{ - 1}}\bm \Sigma {\bm S^{ - 1}}} \right),\]
where $\bm S=\left( {\begin{array}{*{20}{c}}
   {{f_{{\eta _j}}}\left( 0 \right)} & 0  \\
   0 & {{f_{{\eta _j}}}\left( 0 \right){\mu _2}}  \\
\end{array}} \right)$ and $\bm \Sigma=\left( {\begin{array}{*{20}{c}}
   {\tau_j \left( {1 - \tau_j } \right){\nu _0}} & 0  \\
   0 & {\tau_j \left( {1 - \tau_j } \right){\nu _2}}  \\
\end{array}} \right) $.

\noindent \textbf{Proof of Lemma 1.}  To apply the identity (Knight, 1998)
\begin{gather}
{\rho _\tau }\left( {u - v} \right) - {\rho _\tau }\left( u \right) = v\left\{ {I\left( {u \le 0} \right) - \tau } \right\} + \int_0^v {\left\{ {I\left( {u \le t} \right) - I\left( {u \le 0} \right)} \right\}dt}.
 \tag{A.1}
\end{gather}
Minimizing expression (1) is equivalent to minimizing
\[{L_n}\left( \bm\theta  \right) = \sum\limits_{i = 1}^n {\left[ {{\rho _{\tau_j} }\left\{ {{\eta _{ij}} + {d_{ij}} - {\Delta _{ij}}} \right\}{K_{ij}} - {\rho _{\tau_j} }\left\{ {{\eta _{ij}} + {d_{ij}}} \right\}{K_{ij}}} \right]}. \]
Using identity (A.1) and with some straightforward calculations, it follows that
\begin{gather}
\begin{array}{l}
 {L_n}\left( \bm\theta  \right) =  {\frac{1}{{\sqrt {n{h_j}} }}\sum\limits_{i = 1}^n {{K_{ij}}S_{ij}^*u} }  + \frac{1}{{\sqrt {n{h_j}} }}\sum\limits_{i = 1}^n {{K_{ij}}S_{ij}^*{x_{ij}}v}  + {B_{nj}}\left( \bm \theta  \right) \\
  ~~~~~~~~= \bm W_n^{*T}\bm \theta  + {B_{nj}}\left( \bm \theta  \right), \\
 \end{array}
  \tag{A.2}
\end{gather}
where ${B_{nj}}\left( \bm\theta  \right) = \sum\limits_{i = 1}^n {\int_0^{{\Delta _{ij}}} {{K_{ij}}\left\{ {I\left( {{\eta _{ij}} + {d_{ij}} \le t} \right) - I\left( {{\eta _{ij}} + {d_{ij}} \le 0} \right)} \right\}dt} } $. Then
\begin{gather}
\begin{array}{l}
 E\left( {{B_{nj}}\left( \bm\theta  \right)\left| {{\mathcal{X}_j}} \right.} \right) = \sum\limits_{i = 1}^n {\int_0^{{\Delta _{ij}}} {{K_{ij}}\left\{ {{F_{{\eta _j}}}\left( {t - {d_{ij}}} \right) - {F_{{\eta _j}}}\left( { - {d_{ij}}} \right)} \right\}dt} }  \\
  ~~~~~~~~~~~~~~~~~~~~= \sum\limits_{i = 1}^n {\int_0^{{\Delta _{ij}}} {{K_{ij}}\left\{ {{f_{{\eta _j}}}\left( { - {d_{ij}}} \right)t\left( {1 + o\left( 1 \right)} \right)} \right\}dt} }  \\
   ~~~~~~~~~~~~~~~~~~~~= \frac{1}{2}{f_{{\eta _j}}}\left( 0 \right)\sum\limits_{i = 1}^n {{K_{ij}}\Delta _{ij}^2}  + {o_p}\left( 1 \right) \\
   ~~~~~~~~~~~~~~~~~~~~= \frac{1}{2}{\bm\theta ^T}{\bm S_n}\bm\theta  + {o_p}\left( 1 \right), \\
 \end{array}
  \tag{A.3}
\end{gather}
where
\[{\bm S_n} = \left( {\begin{array}{*{20}{c}}
   {\frac{1}{{n{h_j}}}{f_{{\eta _j}}}\left( 0 \right)\sum\limits_{i = 1}^n {{K_{ij}}} } & {\frac{1}{{n{h_j}}}{f_{{\eta _j}}}\left( 0 \right)\sum\limits_{i = 1}^n {{x_{ij}}{K_{ij}}} }  \\
   {\frac{1}{{n{h_j}}}{f_{{\eta _j}}}\left( 0 \right)\sum\limits_{i = 1}^n {{x_{ij}}{K_{ij}}} } & {\frac{1}{{n{h_j}}}{f_{{\eta _j}}}\left( 0 \right)\sum\limits_{i = 1}^n {x_{ij}^2{K_{ij}}} }  \\
\end{array}} \right),\]
and
\begin{gather}
\begin{array}{l}
 Var\left( {{B_{nj}}\left( \bm\theta  \right)\left| {{\mathcal{X}_j}} \right.} \right)\\
  = \sum\limits_{i = 1}^n {Var\left\{ {\left( {\int_0^{{\Delta _{ij}}} {{K_{ij}}\left\{ {I\left( {{\eta _{ij}} + {d_{ij}} \le t} \right) - I\left( {{\eta _{ij}} + {d_{ij}} \le 0} \right)} \right\}dt} } \right)\left| {{\mathcal{X}_j}} \right.} \right\}}  \\
  \le \sum\limits_{i = 1}^n {E\left\{ {{{\left( {\int_0^{{\Delta _{ij}}} {{K_{ij}}\left\{ {I\left( {{\eta _{ij}} + {d_{ij}} \le t} \right) - I\left( {{\eta _{ij}} + {d_{ij}} \le 0} \right)} \right\}dt} } \right)}^2}\left| {{\mathcal{X}_j}} \right.} \right\}}  \\
  \le \sum\limits_{i = 1}^n {K_{ij}^2} \int_0^{\left| {{\Delta _{ij}}} \right|} {\int_0^{\left| {{\Delta _{ij}}} \right|} {\left\{ {{F_{{\eta _j}}}\left( {\left| {{\Delta _{ij}}} \right| - {d_{ij}}} \right) - {F_{{\eta _j}}}\left( { - {d_{ij}}} \right)} \right\}d{t_1}d{t_2}} }  \\
  = o\left( {\sum\limits_{i = 1}^n {K_{ij}^2\Delta _{ij}^2} } \right) = {o_p}\left( 1 \right). \\
 \end{array}
  \tag{A.4}
\end{gather}
Similar to Parzen (1962), we have
\[\frac{1}{{n{h_j}}}\sum\limits_{i = 1}^n {{K_{ij}}x_{ij}^k} \mathop  \to \limits^p {f_j}\left( x \right){\mu _k},\]
where $\mathop  \to \limits^p$ stands for convergence in probability. Thus
\[{\bm S_n}\mathop  \to \limits^p  {f_j}\left( x \right)\left( {\begin{array}{*{20}{c}}
   {{f_{{\eta _j}}}\left( 0 \right)} & 0  \\
   0 & {{f_{{\eta _j}}}\left( 0 \right){\mu _2}}  \\
\end{array}} \right) = {f_j}\left( x \right)\bm S. \]
This together with (A.2)--(A.4), leads to
\[{L_n}\left( \bm\theta  \right) = \frac{1}{2}{f_j}\left( x \right){\bm\theta ^T}\bm S\bm\theta  + \bm W_n^{*T}\bm\theta  + {o_p}\left( 1 \right)\]
Since the convex function ${L_n}\left( \bm\theta  \right) - \bm W_n^{*T}\bm\theta $ converges in probability to the convex function $\frac{1}{2}{f_j}\left( x \right){\bm\theta ^T}\bm S\bm \theta $, it follows from the convexity lemma (Pollard, 1991) that, for any compact
set $\Theta$, the quadratic approximation to ${L_n}\left( \bm\theta  \right)$ holds uniformly for $\bm\theta$ in any compact set, which
leads to
\[\bm{\hat \theta}  =  - \frac{1}{{{f_j}\left( x \right)}}{\bm S^{ - 1}}\bm W_n^* + {o_p}\left( 1 \right).\]
Denote ${S_{ij}} = I\left( {{\eta _{ij}} \le 0} \right) - \tau_j $ and ${\bm W_n} = {\left( {{w_1},{w_2}} \right)^T}$ with ${w_1} = \frac{1}{{\sqrt {n{h_j}} }}\sum\limits_{i = 1}^n {{K_{ij}}{S_{ij}}} $ and ${w_2} = \frac{1}{{\sqrt {n{h_j}} }}\sum\limits_{i = 1}^n {{K_{ij}}{S_{ij}}{x_{ij}}} $.
By the Cram\'{e}r-Wold theorem and central limit theorem, we have
\begin{gather}
\frac{{{\bm W_n}\left| {{\mathcal{X}_j}} \right. - E\left( {{\bm W_n}\left| {{\mathcal{X}_j}} \right.} \right)}}{{\sqrt {Var\left( {{\bm W_n}\left| {{\mathcal{X}_j}} \right.} \right)} }}\mathop  \to \limits^d N\left( {\bm 0,\bm I_{2\times 2}} \right).
  \tag{A.5}
\end{gather}
Note that $E\left( S_{ij} \right) = 0$ and $Var\left( S_{ij} \right) = \tau_j \left( {1 - \tau_j } \right)$. Similar to Parzen (1962), we have $\frac{1}{{ {n{h_j}} }}\sum\limits_{i = 1}^n {K_{ij}^2x_{ij}^k} \mathop  \to \limits^p {f_j}\left( x \right){\nu _k}$. Therefore,
 \[Var\left( {{\bm W_n}\left| {{\mathcal{X}_j}} \right.} \right)\mathop  \to \limits^p {f_j}\left( x \right)\left( {\begin{array}{*{20}{c}}
   {\tau_j \left( {1 - \tau_j } \right){\nu _0}} & 0  \\
   0 & {\tau _j\left( {1 - \tau_j } \right){\nu _2}}  \\
\end{array}} \right)  ={f_j}\left( x \right)\bm \Sigma  .\]
 Combined with (A.5), we have ${\bm W_n}\left| {{\mathcal{X}_j}} \right.\mathop  \to \limits^d N\left( {\bm 0,{f_j}\left( x \right)\bm\Sigma } \right)$. Moreover, we have
\[\begin{array}{l}
 Var\left( {w_1^* - {w_1}\left| {{\mathcal{X}_j}} \right.} \right) = \frac{1}{{n{h_j}}}\sum\limits_{i = 1}^n {K_{ij}^2Var\left( {S_{ij}^* - {S_{ij}}} \right)}  \\
  ~~~~~~~~~~~~~~~~~~~~~~~~\le \frac{1}{{n{h_j}}}\sum\limits_{i = 1}^n {K_{ij}^2\left[ {{F_{{\eta _j}}}\left( {\left| {{d_{ij}}} \right|} \right) - {F_{{\eta _j}}}\left( 0 \right)} \right]}  = {o_p}\left( 1 \right), \\
 \end{array}\]
and
\[\begin{array}{l}
 Var\left( {w_2^* - {w_2}\left| {{\mathcal{X}_j}} \right.} \right) = \frac{1}{{n{h_j}}}\sum\limits_{i = 1}^n {K_{ij}^2x_{ij}^2Var\left( {S_{ij}^* - {S_{ij}}} \right)}  \\
   ~~~~~~~~~~~~~~~~~~~~~~~~\le \frac{1}{{n{h_j}}}\sum\limits_{i = 1}^n {K_{ij}^2x_{ij}^2\left[ {{F_{{\eta _j}}}\left( {\left| {{d_{ij}}} \right|} \right) - {F_{{\eta _j}}}\left( 0 \right)} \right]}  = {o_p}\left( 1 \right) .\\
 \end{array}\]
Thus, we obtain $Var\left( {\bm W_n^* - {\bm W_n}\left| {{\mathcal{X}_j}} \right.} \right) = {o_p}\left( 1 \right)$. So by Slutsky's theorem, conditioning on ${\mathcal{X}_j}$, we have
\[\bm W_n^*\left| {{\mathcal{X}_j}} \right. - E\left( {\bm W_n^*\left| {{\mathcal{X}_j}} \right.} \right)\mathop  \to \limits^d N\left( {\bm 0,{f_j}\left( x \right)\bm \Sigma } \right).\]
Therefore,
\[\bm{\hat \theta}  + \frac{1}{{{f_j}\left( x \right)}}{\bm S^{ - 1}}E\left( {\bm W_n^*\left| {{\mathcal{X}_j}} \right.} \right)\mathop  \to \limits^d N\left( {\bm 0,\frac{1}{{{f_j}\left( x \right)}}{\bm S^{ - 1}}\bm \Sigma {\bm S^{ - 1}}} \right).\]
This completes the proof. $\square$\\

\noindent \textbf{Proof of Theorem 1.}  Let us calculate the conditional bias and variance.
\[\begin{array}{l}
 bias\left( {{{\hat m}_j}\left( x \right)\left| {{\mathcal{X}_j}} \right.} \right) =  - \frac{1}{{\sqrt {n{h_j}} {f_j}\left( x \right)}}\frac{1}{{{f_{{\eta _j}}}\left( 0 \right)}}\frac{1}{{\sqrt {n{h_j}} }}E\left( {\sum\limits_{i = 1}^n {{K_{ij}}S_{ij}^*\left| {{\mathcal{X}_j}} \right.} } \right) \\
   ~~~~~~~~~~~~~~~~~~~~~=  - \frac{1}{{n{h_j}}}\frac{1}{{{f_j}\left( x \right){f_{{\eta _j}}}\left( 0 \right)}}\sum\limits_{i = 1}^n {{K_{ij}}\left[ {{F_{{\eta _j}}}\left( { - {d_{ij}}} \right) - {F_{{\eta _j}}}\left( 0 \right)} \right]}  \\
 ~~~~~~~~~~~~~~~~~~~~~ = \frac{1}{{n{h_j}}}\frac{1}{{{f_j}\left( x \right)}}\sum\limits_{i = 1}^n {{K_{ij}}{d_{ij}}}.  \\
 \end{array}\]
By using the fact that
\[\frac{1}{{n{h_j}}}\sum\limits_{i = 1}^n {{K_{ij}}{d_{ij}}}  = \frac{1}{2}{{f_j}\left( x \right){{\ddot m}_j}\left( x \right)}{\mu _2}h_j^2\left\{ {1 + {o_p}\left( 1 \right)} \right\}.\]
We obtain
\[bias\left( {{{\hat m}_j}\left( x \right)\left| {{\mathcal{X}_j}} \right.} \right) = \frac{1}{2}{{\ddot m}_j}\left( x \right){\mu _2}h_j^2 + {o_p}\left( {h_j^2} \right).\]
Furthermore, the conditional variance of ${{\hat m}_j}\left( x \right)$ is
\[Var\left( {{{\hat m}_j}\left( x \right)\left| {{\mathcal{X}_j}} \right.} \right) = \frac{1}{{n{h_j}{f_j}\left( x \right)}}\frac{{\tau_j \left( {1 - \tau_j } \right){\nu _0}}}{{{{\left( {{f_{{\eta _j}}}\left( 0 \right)} \right)}^2}}} + {o_p}\left( {\frac{1}{{n{h_j}}}} \right).\] $\square$\\

\noindent \textbf{Proof of Theorem 2.} Let ${\alpha _n} = \sqrt {{p_n}} \left( {{n^{{{ - 1} \mathord{\left/
 {\vphantom {{ - 1} 2}} \right.
 \kern-\nulldelimiterspace} 2}}} + {a_n}} \right)$, ${\bm{{\rm{w}}}_o}  = {\left( {{w_{o0}},\bm w_o} \right)^T}$, $\bm w_o=({w_{o1}},...,{w_{o{p_n}}})^T$,  $v ={{{\alpha _n^{-1}}}}\left( {{{\hat w }_0} - {w_{o0}}} \right)$, $ \bm u = \alpha _n^{-1}\left( {\bm{\hat w} - \bm w_o } \right)$ and ${\mathcal{C}_n} = \left\{ {\left( {v,\bm u} \right):{{\left\| {{{\left( {v,{\bm u^T}} \right)}^T}} \right\|}} = C} \right\}$. In order to prove the convergence rate in Theorem 2(i), our aim is to show that for
 any given $\epsilon$ there is a large constant $C$ such that, for a large $n$, we have
 \begin{gather}
P\left\{ {\mathop {\inf }\limits_{(v,\bm u)\in {\mathcal{C}_n}}\mathcal{Q}_n \left({\bm{{\rm{w}}}_o} +\alpha_n(v,\bm u^T)^T \right) \ge {\mathcal{Q}_n }\left( {\bm{{\rm{w}}}_o} \right)} \right\} \ge 1 -\epsilon.
\tag{A.6}
\end{gather}
Using $p_\lambda(0)=0$ and identity (A.1) together with some straightforward calculations, it follows that
\begin{gather}
\begin{array}{l}
{L_n} = \mathcal{Q}_n \left({\bm{{\rm{w}}}_o} +\alpha_n(v,\bm u^T)^T \right) - {\mathcal{Q}_n }\left( {\bm{{\rm{w}}}_o} \right)\\
 ~~~~ \ge \sum\limits_{i = 1}^n {{\alpha _n}\left( {v + {{\hat {\mathcal{M}}}_i}^T\bm u} \right)\left[ {I\left\{ {{Y_i} - {w_{o0}} - {{\hat {\mathcal{M}}}_i}^T{\bm w_o} \le 0} \right\} - \tau } \right]}  \\
  ~~~~~~ + \sum\limits_{i = 1}^n {\int_0^{{\alpha _n}\left( {v + {{\hat {\mathcal{M}}}_i}^T\bm u} \right)} {\left[ {I\left\{ {{Y_i} - {w_{o0}} - {{\hat {\mathcal{M}}}_i}^T{\bm w_o} \le t} \right\} - I\left\{ {{Y_i} - {w_{o0}} - {{\hat {\mathcal{M}}}_i}^T{\bm w_o} \le 0} \right\}} \right]} } dt \\
  ~~~~~~ + n\sum\limits_{j = 1}^{{s_n}} {{p_\lambda }\left( {\left| {{w_{oj}} + {\alpha _n}{u_j}} \right|} \right)}  - n\sum\limits_{j = 1}^{{s_n}} {{p_\lambda }\left( {\left| {{w_{oj}} } \right|} \right)}  \\
  ~~~~ \buildrel \Delta \over = \sqrt n {\alpha _n}(\zeta v + \bm z_n^T\bm u )+ {B_n} + {{\rm P}_n}, \\
\end{array}
  \tag{A.7}
\end{gather}
where
\[\begin{array}{l}
 \zeta  = {n^{{{ - 1} \mathord{\left/
 {\vphantom {{ - 1} 2}} \right.
 \kern-\nulldelimiterspace} 2}}}\sum\limits_{i = 1}^n {\left[ {I\left\{ {{\eta _i} \le {\Delta _i}} \right\} - \tau } \right]},  \\
 \bm z_n = {n^{{{ - 1} \mathord{\left/
 {\vphantom {{ - 1} 2}} \right.
 \kern-\nulldelimiterspace} 2}}}\sum\limits_{i = 1}^n {{{{\hat {\mathcal{M}}}_i}}\left[ {I\left\{ {{\eta _i} \le {\Delta _i}} \right\} - \tau } \right]},  \\
 {B_n} = \sum\limits_{i = 1}^n {\int_0^{{\alpha _n}\left( {v + {{\hat {\mathcal{M}}}_i}^T\bm u} \right)} {\left[ {I\left\{ {{\eta _i} \le t  + {\Delta _i}} \right\} - I\left\{ {{\eta_i} \le  {\Delta _i}} \right\}} \right]} } dt ,\\
{{\rm P}_n} = n\sum\limits_{j = 1}^{{s_n}} {{p_\lambda }\left( {\left| {{w_{oj}} + {\alpha _n}{u_j}} \right|} \right)}  - n\sum\limits_{j = 1}^{{s_n}} {{p_\lambda }\left( {\left| {{w_{oj}}} \right|} \right)},
 \end{array}\]
where ${\Delta _i} = {\left( {{{\hat {\mathcal{M}}}_i} - {\mathcal{M}_i}} \right)^T}\bm w_o  $. Note that, by (C7), $E(\bm z_n^T\bm u)=0$ and $E\left( {{{\left( {\bm z_n^T\bm u} \right)}^2}} \right) = {\bm u^T}E\left( {{\bm z_n}\bm z_n^T} \right)\bm u = C{\bm u^T}{\bm \Psi _n}\bm u \leq C{\bm u^T}\lambda_{max}\left({\bm \Psi _n} \right)\bm u = O\left( {\left\| \bm u \right\|^2} \right)$. Hence, $\bm z_n^T\bm u=O\left( {\left\| \bm u \right\|} \right)$. This combined with (A.7) leads to
\begin{gather}
{L_n}= {B_n} + {{\rm P}_n}+{o_p}\left( {n\alpha _n^2} \right)\left\| \bm u \right\|.
 \tag{A.8}
\end{gather}
By the definition of ${{\hat {\mathcal{M}}}_i}$ and ${{{\mathcal{M}}}_i}$ in section 3 and the consistency result in Theorem 1, we have, for $x_j$ and $j=1,...,p_n$, ${{\hat m}_j}\left( {{x_j}} \right) - {m_j}\left( {{x_j}} \right) = O\left( {h_j^2}+{\left( {n{h_j}} \right)^{{{ - 1} \mathord{\left/
 {\vphantom {{ - 1} 2}} \right.
 \kern-\nulldelimiterspace} 2}}}\right)$. Observe that
\[{{\hat {\mathcal{M}}}_i}{{\hat {\mathcal{M}}}_i}^T - {{ {\mathcal{M}}}_i}{{ {\mathcal{M}}}_i}^T = \left( {{{\hat {\mathcal{M}}}_i} - {{\mathcal{M}}_i}} \right){\mathcal{M}}_i^T + {{\mathcal{M}}_i}{\left( {{{\hat {\mathcal{M}}}_i} - {{\mathcal{M}}_i}} \right)^T} + \left( {{{\hat {\mathcal{M}}}_i} - {{\mathcal{M}}_i}} \right){\left( {{{\hat {\mathcal{M}}}_i} - {{\mathcal{M}}_i}} \right)^T}.\]
Then, for any $\delta$, by Chebyshev's inequality and following the proof of Lemma 8 in Fan and Peng (2004), we have
\[P\left\{ {{{\left\| {\frac{1}{n}\sum\limits_{i = 1}^n {{ {\mathcal{M}}_i} {\mathcal{M}}_i^T - {\bm \Lambda _n}} } \right\|}_F} > \delta } \right\} \le \frac{1}{{{n^2}{\delta ^2}}}E\left[ {\left\| {\sum\limits_{i = 1}^n {{ {\mathcal{M}}_i} {\mathcal{M}}_i^T - n{\bm \Lambda _n}} } \right\|_F^2} \right] = O\left( {{{p_n^2} \mathord{\left/
 {\vphantom {{p_n^2} n}} \right.
 \kern-\nulldelimiterspace} n}} \right) = o\left( 1 \right),\]
 where ${\left\| \bm A \right\|_F}$ stands for Frobenius norm for a real matrix $\bm A = \left( {{a_{ij}}} \right)_{i,j}^{m,n}$, namely, ${\left\|\bm A \right\|_F}={\left( {\sum\limits_{i = 1}^m {\sum\limits_{j = 1}^n {a_{ij}^2} } } \right)^{{1 \mathord{\left/
 {\vphantom {1 2}} \right.
 \kern-\nulldelimiterspace} 2}}}$. Hence, we have
\begin{gather}
{\left\| {\frac{1}{n}\sum\limits_{i = 1}^n {{ {\mathcal{M}}_i} {\mathcal{M}}_i^T - {\bm \Lambda _n}} } \right\|_F} = o_p\left( 1 \right).
\tag{A.9}
\end{gather}
Equation (A.9) and condition (C6) imply that ${\bm u^T}\left(\sum\limits_{i = 1}^n {{\mathcal{M}_i}\mathcal{M}_i^T}/n \right)\bm u $ is asymptotically dominated by ${\bm u^T}\bm\Lambda_n \bm u$. As $np_n\left( {h_j^2 + {{\left( {n{h_j}} \right)}^{{{ - 1} \mathord{\left/
 {\vphantom {{ - 1} 2}} \right.
 \kern-\nulldelimiterspace} 2}}}} \right) = o\left( 1 \right)$ for $j=1,...,p_n$, we can easily prove that
\[{\left\| \frac{1}{n}{\sum\limits_{i = 1}^n {\left( {{{\hat  {\mathcal{M}}}_i} - {\mathcal{M}_i}} \right) {\mathcal{M}}_i^T} } \right\|_F} \to 0,{\left\|\frac{1}{n} {\sum\limits_{i = 1}^n {{ {\mathcal{M}}_i}{{\left( {{{\hat  {\mathcal{M}}}_i} - { {\mathcal{M}}_i}} \right)}^T}} } \right\|_F} \to 0,\]
and
\[{\left\| \frac{1}{n}{\sum\limits_{i = 1}^n {\left( {{{\hat  {\mathcal{M}}}_i} - { {\mathcal{M}}_i}} \right){{\left( {{{\hat  {\mathcal{M}}}_i} - { {\mathcal{M}}_i}} \right)}^T}} } \right\|_F} \to 0.\]
Thus $\frac{1}{n}\sum\limits_{i = 1}^n{{\hat {\mathcal{M}}}_i}{{\hat {\mathcal{M}}}_i}^T - \frac{1}{n}\sum\limits_{i = 1}^n{{ {\mathcal{M}}}_i}{{ {\mathcal{M}}}_i}^T=o_p(1)$.
 Let $\mathbb{M}=({{\hat {\mathcal{M}}}_1^T},...,{{\hat {\mathcal{M}}}_n^T})^T$. Hence, we have
\begin{gather}
\begin{array}{l}
 E\left( {{B_n}}\left| \mathbb{M} \right. \right) = \sum\limits_{i = 1}^n {\int_0^{{{\alpha _n}\left( {v + {{\hat {\mathcal{M}}}_i}^T\bm u} \right)}} {\left[ {{F_\eta }\left\{ {t +  {\Delta _i}} \right\} - {F_\eta }\left\{ { {\Delta _i}} \right\}} \right]} } dt \\
  ~~~~~~~~~~~~~~= \sum\limits_{i = 1}^n {\int_0^{{{\alpha _n}\left( {v + {{\hat {\mathcal{M}}}_i}^T\bm u} \right)}} {\left[ {{f_\eta }\left( { {\Delta _i}} \right)t\left( {1 + o\left( 1 \right)} \right)} \right]} } dt \\
    ~~~~~~~~~~~~~~= \frac{1}{2}\sum\limits_{i = 1}^n {{f_\eta }\left( {0} \right){{\left( {{\alpha _n}\left( {v + {{\hat {\mathcal{M}}}_i}^T\bm u} \right)} \right)}^2}\left( {1 + o\left( 1 \right)} \right)}  \\
  ~~~~~~~~~~~~~~=\frac{1}{2}{f_\eta }\left( 0 \right)n\alpha _n^2\left( {{v^2} + {\bm u^T}{\bm \Psi _n}\bm u + 2v\bm \mu _n^T\bm u} \right)\left( {1 + o_p\left( 1 \right)} \right).\\
 \end{array}
  \tag{A.10}
\end{gather}
Because ${{{p_n^2} \mathord{\left/
 {\vphantom {{p_n^2} n}} \right.
 \kern-\nulldelimiterspace} n}}\rightarrow 0 $, we have ${{{p_n}\log n} \mathord{\left/
 {\vphantom {{{p_n}\log n} n}} \right.
 \kern-\nulldelimiterspace} n} \to 0,\sqrt {{p_n}\log n} {a_n} \to 0$ as $n\rightarrow \infty$, with condition (C5) we
have ${\max _{1 \le i \le n}}\left| {{\alpha _n}\left( {v + {{\hat  {\mathcal{M}}}_i}^T\bm u} \right)} \right| \to 0$. By the Schwarz inequality, it is
not difficult to show that
\begin{gather}
\begin{array}{l}
 Var\left( {{B_n}} \left| \mathbb{M} \right.\right) \\
 \le \sum\limits_{i = 1}^n {E{{\left( {\int_0^{{\alpha _n}\left( {v + {{\hat {\mathcal{M}}}_i}^T\bm u} \right)} {\left[ {I\left\{ {{\eta _i} \le t + {\Delta _i}} \right\} - I\left\{ {{\eta _i} \le {\Delta _i}} \right\}} \right]dt\left| \mathbb{M} \right.} } \right)}^2}}  \\
  \le \sum\limits_{i = 1}^n {{\alpha _n}\left| {v + {{\hat {\mathcal{M}}}_i}^T\bm u} \right|\left| {\int_0^{{\alpha _n}\left( {v + {{\hat {\mathcal{M}}}_i}^T\bm u} \right)} {E\left[ {I\left\{ {{\eta _i} \le t + {\Delta _i}} \right\} - I\left\{ {{\eta _i} \le {\Delta _i}} \right\}} \right]^2dt\left|\mathbb{ M} \right.} } \right|}  \\
  = o\left( 1 \right)\sum\limits_{i = 1}^n {{{\left( {{\alpha _n}\left( {v + {{\hat{\mathcal{M}}}_i}^T\bm u} \right)} \right)^2}}}  \\
  = {o_p}\left(p_n \right) .\\
 \end{array}
  \tag{A.11}
\end{gather}
From (A.9)--(A.11), we have
\begin{gather}{L_n} = \frac{1}{2}{f_\eta }\left( 0 \right)n\alpha _n^2\left( {{v^2} + {\bm u^T}{\bm\Psi _n}\bm u + 2v\bm \mu _n^T\bm u} \right)+{o_p}\left( {n\alpha _n^2} \right)\left\| \bm u \right\|+{{\rm P}_n}.
 \tag{A.12}
\end{gather}
For ${{\rm P}_n}$, by the condition (C9) and Taylor's expansion for the penalty function, we have
\begin{gather}
\begin{array}{l}
 {{\rm P}_n} = n\sum\limits_{j = 1}^{{s_n}} {\left[ {{p_\lambda }\left( {\left| {{w_{oj}} + {\alpha _n}{u_j}} \right|} \right) - {p_\lambda }\left( {\left| {{w_{oj}}} \right|} \right)} \right]}  \\
  ~~~= n\sum\limits_{j = 1}^{{s_n}} {{\dot{p}_\lambda }\left( {\left| {{w_{oj}}} \right|} \right){\mathop{\rm sgn}} \left( {{w_{oj}}} \right){\alpha _n}{u_j}}  + \frac{1}{2}n\sum\limits_{j = 1}^{{s_n}} {{{\ddot{p}}_\lambda }\left( {\left| {{w_{oj}}} \right|} \right){{\left( {{\alpha _n}{u_j}} \right)}^2}\left\{ {1 + o\left( 1 \right)} \right\}}  \\
   ~~~\le n{\alpha _n}\sqrt {{s_n}} {\max _{1 \le j \le {s_n}}}\left| {{{\dot{p}}_\lambda }\left( {\left| {{w_{oj}}} \right|} \right)} \right|\left\| \bm u \right\| + n\alpha _n^2{\max _{1 \le j \le {s_n}}}\left| {{{\ddot{p}}_\lambda }\left( {\left| {{w_{oj}}} \right|} \right)} \right|{\left\| \bm u \right\|^2}\left\{ {1 + o\left( 1 \right)} \right\} \\
   ~~~= O\left( {n{\alpha _n}{a_n}\sqrt {{s_n}} } \right)\left\| \bm u \right\| + O\left( {n\alpha _n^2{b_n}} \right){\left\| \bm u \right\|^2} \\
  ~~~ = O\left( {n\alpha _n^2} \right){\left\| \bm u \right\|}+o\left( {n\alpha _n^2} \right){\left\| \bm u \right\|^2}. \\
 \end{array}
  \tag{A.13}
\end{gather}
It follows from (A.12) and (A.13) that $L_n$ in (A.7) is dominated by the positive quadratic term $\frac{1}{2}{f_\eta }\left( 0 \right)n\alpha _n^2\left( {{v^2} + {\bm u^T}{\bm \Psi _n}\bm u + 2v\bm\mu _n^T\bm u} \right)$ when
a sufficiently large $C$ is chosen. Therefore, (A.6) holds and this completes the proof of Theorem 2 (i).

(ii) Let ${\bm{{\rm{\hat w}}}_n(1)}$ and ${\bm{{\rm{\hat w}}}_n(2)}$ be the estimators of ${\bm{{\rm{w}}}_n(1)}$ and ${\bm{{\rm{w}}}_n(2)}$, respectively, where${\bm{{\rm{w}}}_n(1)}=(w_{1},...,w_{s_n})^T$ and ${\bm{{\rm{w}}}_n(2)}=(w_{s_n+1},...,w_{p_n})^T$. To prove Theorem 2(ii), it suffices to show that for any constant $C$ and any given $(w_0,{\bm{{\rm{w}}}_n^T(1)})^T$ satisfying $\left\|(w_0,{\bm{{\rm{w}}}_n^T(1)})^T-(w_0,{\bm{{\rm{w}}}_o^T(1))^T} \right\|=O_p\left({\alpha _n}\right) $, where ${\alpha _n} = \sqrt {{p_n}} \left( {{n^{{{ - 1} \mathord{\left/
 {\vphantom {{ - 1} 2}} \right.
 \kern-\nulldelimiterspace} 2}}} + {a_n}} \right)$, we have
\begin{gather}
{\mathcal{Q}_n}\left([ w_{0},{  \bm{{\rm{w}}}_n^T\left( 1 \right),{\bm 0^T}} ]^T\right) = \mathop {\min }\limits_{\left\| {{\bm{{\rm{w}}}_n}\left( 2 \right)} \right\| \le C{\alpha _n}} {{\mathcal{Q}}_n}\left( [w_0,{\bm{{\rm{w}}}_n^T\left( 1 \right),\bm{{\rm{w}}}_n^T\left( 2 \right)} ]^T\right).
  \tag{A.14}
\end{gather}
 By (A.14),  it is easy to prove that ${\bm{{\rm{\hat w}}}_n(2)}=0$.

To prove (A.14), it is sufficient to show that, with probability approaching
one, for any $(p_n+1)$-dimensional vector ${\bm{{\rm{ w}}}_n}=(w_0,{\bm{{\rm{ w}}}_n^T(1)},{\bm{{\rm{w}}}_n^T(2)})^T$ with $(w_0,{\bm{{\rm{w}}}_n^T(1)})^T$ satisfying $\left\|(w_0,{\bm{{\rm{w}}}_n^T(1)})^T-(w_0,{\bm{{\rm{w}}}_o^T(1))^T} \right\|=O_p\left({\alpha _n}\right) $
and for some small $\epsilon_n=C{\alpha _n}$ and $j=s_n+1,...,p_n$,
\begin{gather}
\frac{{\partial \mathcal{Q}_n\left( {\bm{{\rm{w}}}_{n}}\right)}}{{\partial {w_j}}} > 0,~~0 < {w_j} < {\epsilon _n},
  \tag{A.15}
\end{gather}
and
\begin{gather}
\frac{{\partial \mathcal{Q}_n\left( {\bm{{\rm{w}}}_{n}}\right)}}{{\partial {w_j}}} < 0,~~-{\epsilon _n} < {w_j} < 0,
  \tag{A.16}
\end{gather}
 Taking the first derivative of $\mathcal{Q}_n\left( {\bm{{\rm{w}}}_{n}}\right)$ at any differentiable point ${\bm{{\rm{ w}}}_n}=(w_0,w_1,...,w_{p_n})^T$ with respect to ${{{{ w}}}_{j}}, j=s_n+1,...,p_n$, we have
\[\frac{{\partial \mathcal{Q}_n\left( {\bm{{\rm{w}}}_{n}}\right)}}{{\partial {w_j}}}=  - \sum\limits_{i = 1}^n {{\psi _\tau }\left\{ {{Y_i} - {w_0} - \sum\limits_{j = 1}^{{p_n}} {{{\hat m}_j}\left( {{X_{ij}}} \right){w_j}} } \right\}{{\hat m}_j}\left( {{X_{ij}}} \right)}  + n{{\dot{p}}_\lambda }\left( {\left| {{w_j}} \right|} \right){\mathop{\rm sgn}} \left( {{w_j}} \right)\]
for $j=1,...,p_n$, where ${\psi _\tau }\left( u \right)=\tau  - I\left( {u < 0} \right)$ and
\[\begin{array}{l}
  - \sum\limits_{i = 1}^n {{\psi _\tau }\left\{ {{Y_i} - {w_0} - \sum\limits_{j = 1}^{{p_n}} {{{\hat m}_j}\left( {{X_{ij}}} \right){w_j}} } \right\}{{\hat m}_j}\left( {{X_{ij}}} \right)}  \\
  =  - \sum\limits_{i = 1}^n {{\psi _\tau }\left\{ {{\eta _i} - \left( {{w_0} - {w_{o0}}} \right) - {\varsigma _i}} \right\}{{\hat m}_j}\left( {{X_{ij}}} \right)}  \\
  =  \sum\limits_{i = 1}^n {\left\{ {I\left( {{\eta _i} \le 0} \right) - \tau } \right\}{{\hat m}_j}\left( {{X_{ij}}} \right)}  + \sum\limits_{i = 1}^n {\left[ {I\left\{ {{\eta _i} \le \left( {{w_0} - {w_{o0}}} \right) + {\varsigma_i}} \right\} - \left\{ {I\left( {{\eta _i} \le 0} \right)} \right\}} \right]{{\hat m}_j}\left( {{X_{ij}}} \right)}  \\
  \buildrel \Delta \over = I + II, \\
 \end{array}\]
where ${\varsigma _i}= \sum\limits_{j = 1}^{{p_n}} {\left\{ {\left[ {{{\hat m}_j}\left( {{X_{ij}}} \right) - {m_j}\left( {{X_{ij}}} \right)} \right]{w_{oj}} + {{\hat m}_j}\left( {{X_{ij}}} \right)\left[ {{w_j} - {w_{oj}}} \right]} \right\}} $.
As in the proof of Theorem 2 (i) and Theorem 1, it is easy to prove that
\[I=O_p\left( \sqrt {n{p_n}}  \right)~{\rm{and}}~II=O_p\left( \sqrt {n{p_n}}  \right).\]
Hence, we have
\begin{gather}
\begin{array}{l}
 \frac{{\partial \mathcal{Q}_n\left( {\bm{{\rm{w}}}_{n}}\right)}}{{\partial {w_j}}}= {O_p}\left( {\sqrt {n{p_n}} } \right) + n{{\dot{p}}_\lambda }\left( {\left| {{w_j}} \right|} \right){\mathop{\rm sgn}} \left( {{w_j}} \right) \\
  ~~~~~~~~~~= n\lambda \left\{ {{O_p}\left( {{{\sqrt {{{{p_n}} \mathord{\left/
 {\vphantom {{{p_n}} n}} \right.
 \kern-\nulldelimiterspace} n}} } \mathord{\left/
 {\vphantom {{\sqrt {{{{p_n}} \mathord{\left/
 {\vphantom {{{p_n}} n}} \right.
 \kern-\nulldelimiterspace} n}} } \lambda }} \right.
 \kern-\nulldelimiterspace} \lambda }} \right) + {\lambda ^{ - 1}}{{\dot{p}}_\lambda }\left( {\left| {{w_j}} \right|} \right){\mathop{\rm sgn}} \left( {{w_j}} \right)} \right\} .\\
 \end{array}
   \tag{A.17}
\end{gather}
Whereas $\lim {\inf _{n \to \infty }}\lim {\inf _{\left| {{w_j}} \right| \to {0_ + }}}{\dot{p}_{\lambda}}\left( \left| {{w_j}} \right| \right)/\lambda > 0$ by condition (C8) and ${{\sqrt {{{{p_n}} \mathord{\left/
 {\vphantom {{{p_n}} n}} \right.
 \kern-\nulldelimiterspace} n}} } \mathord{\left/
 {\vphantom {{\sqrt {{{{p_n}} \mathord{\left/
 {\vphantom {{{p_n}} n}} \right.
 \kern-\nulldelimiterspace} n}} } \lambda }} \right.
 \kern-\nulldelimiterspace} \lambda } \to 0$, the sign of the derivative is completely determined by that of $w_j$ .
Hence, we can show that (A.15) and (A.16) hold by using (A.17). This completes the proof of Theorem 2 (ii).

(iii) It can be shown easily that there exists a $(\hat w_0, {\bm{{\rm{\hat w}}}_n^T(1)})^T$ in Theorem 2 that is a
$\sqrt {{n \mathord{\left/
 {\vphantom {n {{p_n}}}} \right.
 \kern-\nulldelimiterspace} {{p_n}}}} $ consistent local minimizer of ${\mathcal{Q}_n}\left((w_0, {  \bm{{\rm{w}}}_n^T\left( 1 \right),{\bm 0^T}} )^T\right)$, and satisfies the equations
\begin{gather}
 n^{-1}\sum\limits_{i = 1}^n {\hat{\mathcal{M}}_{1i}{\psi _\tau }\left\{ {{Y_i} - { \hat w_{0}} -\hat{\mathcal{M}}_{1i}^T\bm{{\rm{\hat w}}}_n\left( 1 \right)} \right\}}   =  {\bm c_n} + {\bm \Gamma _n}\left( \bm{{\rm{\hat w}}}_n\left( 1 \right)-\bm{{\rm{w}}}_{o}\left( 1 \right) \right)\left\{ {1 + {o_p}(1)} \right\},
 \tag{A.18}
\end{gather}
and
\begin{gather}
- n^{-1}\sum\limits_{i = 1}^n {{\psi _\tau }\left\{ {{Y_i} - {{\hat w}_{0}} - \hat{\mathcal{M}}_{1i}^T{\bm{{\rm{\hat w}}}_n}(1)} \right\}}  = 0.
 \tag{A.19}
\end{gather}
where ${{\hat {\mathcal{M}}}_{1i}} = {\left( {{{\hat m}_1}\left( {{X_{i1}}} \right),...,{{\hat m}_{s_n}}\left( {{X_{i{s_n}}}} \right)} \right)^T}$. We can write
\begin{gather}
 n^{-1}\sum\limits_{i = 1}^n {\hat{\mathcal{M}}_{1i}{\psi _\tau }\left\{ {{Y_i} - { \hat w_{0}} -\hat{\mathcal{M}}_{1i}^T\bm{{\rm{\hat w}}}_n\left( 1 \right)} \right\}}
  =  - {n^{{{ - 1} \mathord{\left/
 {\vphantom {{ - 1} 2}} \right.
 \kern-\nulldelimiterspace} 2}}}{\bm z_{n1}} + {B_{n1}} + {B_{n2}},
  \tag{A.20}
\end{gather}
where
\[\begin{array}{l}
 {\bm z_{n1}} = {n^{{{ - 1} \mathord{\left/
 {\vphantom {{ - 1} 2}} \right.
 \kern-\nulldelimiterspace} 2}}}\sum\limits_{i = 1}^n {{\hat{\mathcal{M}}_{1i}}\left\{ {I\left( {{\eta _i} \le 0} \right) - \tau } \right\}},  \\
 {B_{n1}} = {n^{ - 1}}\sum\limits_{i = 1}^n {{{\hat{\mathcal{M}}}_{1i}}\left\{ {{F_\eta }\left( 0 \right) - {F_\eta }\left( {{\zeta _i}} \right)} \right\}} , \\
 {B_{n2}} = {n^{ - 1}}\sum\limits_{i = 1}^n {{\hat{\mathcal{M}}_{1i}}\left\{ {\left[ {I\left( {{\eta _i} \le 0} \right) - I\left( {{\eta _i} \le {\zeta _i}} \right)} \right] - \left[ {{F_\eta }\left( 0 \right) - {F_\eta }\left( {{\zeta _i}} \right)} \right]} \right\}} , \\
 \end{array}\]
 where ${\zeta _i} = \left( {{{\hat w}_0} - {w_{o0}}} \right) + \sum\limits_{j = 1}^{{s_n}} {\left[ {{{\hat m}_j}\left( {{X_{ij}}} \right)\left( {{{\hat w}_j} - {w_{oj}}} \right) + \left( {{{\hat m}_j}\left( {{X_{ij}}} \right) - {m_j}\left( {{X_{ij}}} \right)} \right){w_{oj}}} \right]} $.

Taking Taylor¡¯s explanation for ${F_\eta }\left({\zeta _i} \right)$ at 0 gives
\begin{gather}
\begin{array}{l}
 {B_{n1}} = {n^{ - 1}}\sum\limits_{i = 1}^n {{\hat{\mathcal{M}}_{1i}}\left\{ {{F_\eta }\left( 0 \right) - {F_\eta }\left( {{\zeta _i}} \right)} \right\}}  \\
  ~~~~~=  - {n^{ - 1}}{f_\eta }\left( 0 \right)\sum\limits_{i = 1}^n {{\hat{\mathcal{M}}_{1i}}{\zeta _i}\left\{ {1 + {o_p}(1)} \right\}}  \\
  ~~~~~ =  - {n^{ - 1}}{f_\eta }\left( 0 \right)\left[ {\sum\limits_{i = 1}^n {{\hat{\mathcal{M}}_{1i}}\left( {{{\hat w}_0} - {w_{o0}}} \right)}  + \sum\limits_{i = 1}^n {{\hat{\mathcal{M}}_{1i}}\hat{\mathcal{M}}_{1i}^T\left( {{\bm{{\rm{\hat w}}}_n}(1) - {\bm{{\rm{ w}}}_o}(1)} \right)} } \right. \\
  ~~~~~~~~\left. { + \sum\limits_{i = 1}^n {{\hat{\mathcal{M}}_{1i}}{{\left( {{\hat{\mathcal{M}}_{1i}} - {{\mathcal{M}}_{1i}}} \right)}^T}{\bm{{\rm{ w}}}_o}(1)} } \right]\left\{ {1 + {o_p}(1)} \right\} \\
  ~~~~~=- n^{ - 1}{f_\eta }\left( 0 \right)\left[ { {n\bm\mu_{n1}\left( {{{\hat w}_0} - {w_{o0}}} \right)}  +  {n{\bm\Psi _{n1}}\left( {{\bm{{\rm{\hat w}}}_n}(1) - {\bm{{\rm{ w}}}_o}(1)} \right)} } \right. \\
  ~~~~~~~\left. { + \sum\limits_{i = 1}^n {{\hat{\mathcal{M}}_{1i}}{{\left( {{\hat{\mathcal{M}}_{1i}} - {{\mathcal{M}}_{1i}}} \right)}^T}{\bm{{\rm{ w}}}_o}(1)} } \right]\left\{ {1 + {o_p}(1)} \right\}.\\
 \end{array}
   \tag{A.21}
\end{gather}
By direct calculation of the mean and variance, we can show, as in Jiang et al. (2001), that $B_{n2}=o_p(\alpha_n)$. This combined with (A.20) and (A.21) leads to
\begin{gather}
\begin{array}{l}
 \left[ {{f_\eta }\left( 0 \right) {\bm\Psi _{n1}}  + {\bm\Gamma _n}} \right]\left( {{\bm{{\rm{\hat w}}}_n}(1) - {\bm{{\rm{ w}}}_o}(1)} \right)\left\{ {1 + {o_p}(1)} \right\}+{n^{{{ - 1} \mathord{\left/
 {\vphantom {{ - 1} 2}} \right.
 \kern-\nulldelimiterspace} 2}}}{\bm z_{n1}} + {\bm c_n} \\
  =   - {n^{ - 1}}{f_\eta }\left( 0 \right)\left[ { {n\bm\mu_{n1}\left( {{{\hat w}_0} - {w_{o0}}} \right) + \sum\limits_{i = 1}^n {{\hat{\mathcal{M}}_{1i}}{{\left( {{\hat{\mathcal{M}}_{1i}} - {{\mathcal{M}}_{1i}}} \right)}^T}{\bm{{\rm{ w}}}_o}(1)} } } \right]\left\{ {1 + {o_p}(1)} \right\}. \\
 \end{array}
 \tag{A.22}
\end{gather}
Similarly, (A.19) can be simplified as
\begin{gather}
\begin{array}{l}
 {n^{ - 1}}\sum\limits_{i = 1}^n {\left\{ {\tau  - I\left( {{\eta _i} \le 0} \right)} \right\}}  = {f_\eta }\left( 0 \right)\left( {{{\hat w}_0} - {w_{o0}}} \right) \\
  + {n^{ - 1}}{f_\eta }\left( 0 \right)\sum\limits_{i = 1}^n {\left[ {{{\left( {{\hat{\mathcal{M}}_{1i}} - {{\mathcal{M}}_{1i}}} \right)}^T}{\bm{{\rm{ w}}}_o}\left( 1 \right) + \hat{\mathcal{M}}_{1i}^T\left( {{\bm{{\rm{\hat w}}}_n}\left( 1 \right) - {\bm{{\rm{ w}}}_o}\left( 1 \right)} \right)} \right]\left\{ {1 + {o_p}(1)} \right\}}  \\
 \end{array}
\tag{A.23}
\end{gather}
Solving (A.22) and (A.23), we obtain that
\begin{gather}
\begin{array}{l}
 \left( {{\bm{\hat w}_n}(1) - {\bm w_o}(1)} \right) \\
 =  - {\left[ {{f_\eta }\left( 0 \right)\bm G_{n1} + {\bm\Gamma _n}} \right]^{ - 1}}{n^{ - 1}}{f_\eta }\left( 0 \right)\sum\limits_{i = 1}^n {{\hat{\mathcal{M}}_{1i}}{{\left( {{\hat{\mathcal{M}}_{1i}} - {{\mathcal{M}}_{1i}}} \right)}^T}{\bm{{\rm{ w}}}_o}(1)}  \\
  ~~- {\left[ {{f_\eta }\left( 0 \right)\bm G_{n1} + {\bm\Gamma _n}} \right]^{ - 1}}{n^{{{ - 1} \mathord{\left/
 {\vphantom {{ - 1} 2}} \right.
 \kern-\nulldelimiterspace} 2}}}{\bm z_{n1}} \\
 ~~+ {\left[ {{f_\eta }\left( 0 \right)\bm G_{n1} + {\bm\Gamma _n}} \right]^{ - 1}}{n^{ - 1}}{f_\eta }\left( 0 \right){\bm\mu _{n1}}\sum\limits_{i = 1}^n {{{\left( {{\hat{\mathcal{M}}_{1i}} - {{\mathcal{M}}_{1i}}} \right)}^T}{\bm{{\rm{ w}}}_o}\left( 1 \right)}  \\
 ~~ - {\left[ {{f_\eta }\left( 0 \right)\bm G_{n1} + {\bm\Gamma _n}} \right]^{ - 1}}\left( {{\bm c_n} + {n^{ - 1}}\sum\limits_{i = 1}^n {\left\{ {\tau  - I\left( {{\eta _i} \le 0} \right)} \right\}{\bm\mu _{n1}}} } \right) + {o_p}\left( {{n^{{{ - 1} \mathord{\left/
 {\vphantom {{ - 1} 2}} \right.
 \kern-\nulldelimiterspace} 2}}}} \right) \\
  = {\Pi _{n1}} + {\Pi _{n2}} + {\Pi _{n3}} + {\Pi _{n4}} + {o_p}\left( {{n^{{{ - 1} \mathord{\left/
 {\vphantom {{ - 1} 2}} \right.
 \kern-\nulldelimiterspace} 2}}}} \right) .\\
 \end{array}
  \tag{A.24}
\end{gather}
  By Theorem 1, we have, uniformly for $x_k\in \mathcal{C}_k$,
\begin{gather}
{{\hat m}_k}\left( {{x_k}} \right) - {m_k}\left( {{x_k}} \right) = \frac{1}{{n{h_k}}}f_k^{ - 1}\left( {{x_k}} \right)f_{{\eta _k}}^{ - 1}\left( 0 \right)\sum\limits_{t = 1}^n {K\left( {\frac{{{X_{tk}} - {x_k}}}{{{h_k}}}} \right){\psi _{\tau_k} }\left( {{\eta _{tk}}} \right)} \left\{ {1 + {o_p}\left( 1 \right)} \right\}.
 \tag{A.25}
\end{gather}
Then, we have
\begin{gather}
{\Pi _{n1}}\mathop  \to \limits^p {-\left[ {{f_\eta }\left( 0 \right)\bm G_{n1} + {\bm\Gamma _n}} \right]^{ - 1}}{n^{ - 1}}{f_\eta }\left( 0 \right){\Pi _{n5}},
\tag{A.26}
\end{gather}
where
\[\begin{array}{l}
 {\Pi _{n5}}= \left\{ {\sum\limits_{i = 1}^n {{m_j}\left( {{X_{ij}}} \right)\sum\limits_{k = 1}^{{s_n}} {{w_{ok}}\left[ {{{\hat m}_k}\left( {{X_{ik}}} \right) - {m_k}\left( {{X_{ik}}} \right)} \right]} } } \right\}_{j = 1,...,{s_n}}^T \\
  ~~~~~= \left\{ {\sum\limits_{i = 1}^n {{m_j}\left( {{X_{ij}}} \right)\sum\limits_{k = 1}^{{s_n}} {{w_{ok}}} } } \right. \\
  ~~~~~~~\times \left. {\left[ {\frac{1}{{n{h_k}}}f_k^{ - 1}\left( {{X_{ik}}} \right)f_{{\eta _k}}^{ - 1}\left( 0 \right)\sum\limits_{t = 1}^n {K\left( {\frac{{{X_{tk}} - {X_{ik}}}}{{{h_k}}}} \right){\psi _{\tau_k} }\left( {{\eta _{tk}}} \right)} \left\{ {1 + {o_p}\left( 1 \right)} \right\}} \right]} \right\}_{j = 1,...,{s_n}}^T \\
 ~~~~~ = \left\{ { {\sum\limits_{t = 1}^{n}\sum\limits_{k = 1}^{{s_n}} {{w_{ok}}{\psi _{\tau_k} }\left( {{\eta _{tk}}} \right)} } } \right. \\
  ~~~~~~~\times \left. {\left[ {\frac{1}{{n{h_k}}}\sum\limits_{i = 1}^n {{m_j}\left( {{X_{ij}}} \right)f_k^{ - 1}\left( {{X_{ik}}} \right)f_{{\eta _k}}^{ - 1}\left( 0 \right){K\left( {\frac{{{X_{tk}} - {X_{ik}}}}{{{h_k}}}} \right)} \left\{ {1 + {o_p}\left( 1 \right)} \right\}} } \right]} \right\}_{j = 1,...,{s_n}}^T, \\
 \end{array}\]
where $f_k(\cdot)$ is the marginal density function of $X_{ik}$. If $k=j$, we have
\begin{gather}
\frac{1}{{n{h_j}}}\sum\limits_{i = 1}^n {{m_j}\left( {{X_{ij}}} \right)f_j^{ - 1}\left( {{X_{ij}}} \right)f_{{\eta _j}}^{ - 1}\left( 0 \right)K\left( {\frac{{{X_{tj}} - {X_{ij}}}}{{{h_j}}}} \right)}  = {m_j}\left( {{X_{tj}}} \right) + {o_p}\left( 1 \right).
\tag{A.27}
\end{gather}
If $k\neq j$, we have
\begin{gather}
\frac{1}{{n{h_k}}}\sum\limits_{i = 1}^n {{m_j}\left( {{X_{ij}}} \right)f_k^{ - 1}\left( {{X_{ik}}} \right)f_{{\eta _k}}^{ - 1}\left( 0 \right)K\left( {\frac{{{X_{tk}} - {X_{ik}}}}{{{h_k}}}} \right)}  = {\beta _{jk}}\left( {{X_{tk}}} \right) + {o_p}\left( 1 \right),
\tag{A.28}
\end{gather}
where ${\beta _{jk}}\left( {{X_{tk}}} \right) = E\left( {{m_j}\left( {{X_{tj}}} \right)\left| {{X_{tk}}} \right.} \right)$. Then, by (A.27)
 and (A.28) and noting that ${\beta _{jj}}\left( {{X_{tj}}} \right) = {m_j}\left( {{X_{tj}}} \right)$, we have
\begin{gather}
\begin{array}{l}
{\Pi _{n5}} = \left( {\sum\limits_{t = 1}^n {\sum\limits_{k = 1}^{{s_n}} {{w_{ok}}{\psi _{\tau_k} }\left( {{\eta _{tk}}} \right){\beta _{jk}}\left( {{X_{tk}}} \right)} } } \right)_{j = 1,...,{s_n}}^T + {O_p}\left( {ns_n\mathop {\max }\limits_{1 \le j \le {s_n}} \left( {h_j^2+(nh_j)^{-1/2}} \right)}  \right) \\
  ~~~~~= {\left( {\sum\limits_{t = 1}^n {{{\tilde \eta }_{t1}},...,\sum\limits_{t = 1}^n {{{\tilde \eta }_{t{s_n}}}} } } \right)^T} + {O_p}\left( {ns_n\mathop {\max }\limits_{1 \le j \le {s_n}} \left( {h_j^2+(nh_j)^{-1/2}} \right)} \right), \\
 \end{array}
  \tag{A.29}
\end{gather}
where ${{\tilde \eta }_{tj}} = \sum\limits_{k = 1}^{{s_n}} {{w_{ok}}{\psi _{\tau_k} }\left( {{\eta _{tk}}} \right){\beta _{jk}}\left( {{X_{tk}}} \right)} $. By (A.26) and (A.29), we have
\begin{gather}
\begin{array}{l}
{\Pi _{n1}}={-\left[ {{f_\eta }\left( 0 \right)\bm G_{n1} + {\bm\Gamma _n}} \right]^{ - 1}}{n^{ - 1}}{f_\eta }\left( 0 \right){\left( {\sum\limits_{i = 1}^n {{{\tilde \eta }_{i1}},...,\sum\limits_{i = 1}^n {{{\tilde \eta }_{i{s_n}}}} } } \right)^T} +\\
~~~~~~~~~{O_p}\left( {s_n\mathop {\max }\limits_{1 \le j \le {s_n}} \left( {h_j^2+(nh_j)^{-1/2}} \right)} \right).
 \end{array}
\tag{A.30}
\end{gather}
We next consider ${\Pi _{n2}}$. Observe that
\[\begin{array}{l}
 {\Pi _{n2}} = {\left[ {{f_\eta }\left( 0 \right)\bm G_{n1} + {\bm\Gamma _n}} \right]^{ - 1}}{n^{ - 1}}\sum\limits_{i = 1}^n {{{\hat {\mathcal{M}}}_{1i}}{\psi _\tau }\left( {{\eta _i}} \right)}  \\
 ~~~~~ = {\left[ {{f_\eta }\left( 0 \right)\bm G_{n1}+ {\bm\Gamma _n}} \right]^{ - 1}}{n^{ - 1}}\\
  ~~~~~~~~~\times\left[ {\sum\limits_{i = 1}^n {{{\mathcal{M}}_{1i}}{\psi _\tau }\left( {{\eta _i}} \right) + \sum\limits_{i = 1}^n {\left( {\hat{{\mathcal{M}}_{1i}} - {{\mathcal{M}}_{1i}}} \right){\psi _\tau }\left( {{\eta _i}} \right)} } } \right]\left\{ {1 + {o_p}\left( 1 \right)} \right\}. \\
 \end{array}\]
We can show that the leading term of $ {\Pi _{n2}}$ is ${\left[ {{f_\eta }\left( 0 \right)\bm G_{n1} + {\bm\Gamma _n}}  \right]^{ - 1}}{n^{ - 1}}\sum\limits_{i = 1}^n {{{\mathcal{M}}_{1i}}{\psi _\tau }\left( {{\eta _i}} \right)} $.
\begin{gather}
 {\Pi _{n2}} = {\left[ {{f_\eta }\left( 0 \right)\bm G_{n1} + {\bm\Gamma _n}} \right]^{ - 1}}{n^{ - 1}}{\left( {\sum\limits_{i = 1}^n {{{\bar \eta }_{i1}},...,\sum\limits_{i = 1}^n {{{\bar \eta }_{i{s_n}}}} } } \right)^T}\left\{ {1 + {o_p}\left( 1 \right)} \right\},
\tag{A.31}
\end{gather}
and ${\Pi _{n3}} = {o_p}\left( {{n^{{{ - 1} \mathord{\left/
 {\vphantom {{ - 1} 2}} \right.
 \kern-\nulldelimiterspace} 2}}}} \right)$. This combined with (A.30) and (A.31) leads to
\[\begin{array}{l}
\sqrt {n}\left[ {{f_\eta }\left( 0 \right)\bm G_{n1} + {\bm\Gamma _n}} \right]\left( {{\bm{{\rm{\hat w}}}_n}(1) -
{\bm{{\rm{w}}}_o}(1)} \right) +\sqrt {n}\bm c_n=
{n^{ - {1 \mathord{\left/
 {\vphantom {1 2}} \right.
 \kern-\nulldelimiterspace} 2}}}\sum\limits_{i = 1}^n {\left( {{\bm\xi _i}{\psi _\tau }\left( {{\eta _i}} \right) - {f_\eta }\left( 0 \right){\bm{\tilde \eta }_i}} \right)} .
 \end{array}\]
where $\bm e_n^T\bm\Upsilon _{n1}^{-1/2}  {n^{ - {1 \mathord{\left/
 {\vphantom {1 2}} \right.
 \kern-\nulldelimiterspace} 2}}}\sum\limits_{i = 1}^n {\left( {{\bm\xi _i}{\psi _\tau }\left( {{\eta _i}} \right) - {f_\eta }\left( 0 \right){\bm{\tilde \eta }_i}} \right)}\mathop  \to \limits^d N\left( 0,1 \right).$
 It follows that
\[\begin{array}{l}
\sqrt n \bm e_n^T\bm\Upsilon _{n1}^{{{ - 1} \mathord{\left/
 {\vphantom {{ - 1} 2}} \right.
 \kern-\nulldelimiterspace} 2}}{\left(  {{f_\eta }\left( 0 \right)\bm G_{n1} + {\Gamma _n}} \right)}\left[ {\left( {{\bm{{\rm{\hat w}}}_n}(1) - {\bm{{\rm{ w}}}_o}(1)} \right) + {{\left(  {{f_\eta }\left( 0 \right)\bm G_{n1} + {\bm\Gamma _n}}\right)}^{ - 1}}{\bm c_n}} \right] \mathop  \to \limits^d N\left( 0,1\right).\\
  \end{array}\] $\square$\\

\noindent  \textbf{References}

\vspace{1mm}
\noindent Ando, T., Li, K.-C., 2014. A model-averaging approach for high-dimensional regression.
Journal of the American Statistical Association, 109, 254--265.
\vspace{1mm}

\vspace{1mm}
\noindent  Chen, J., Li, D., Linton, O., Lu, Z., 2016. Semiparametric dynamic portfolio choice with multiple conditioning
variables. Journal of Econometrics, 194, 309--318.
\vspace{1mm}

\vspace{1mm}
\noindent  Chen, J., Li, D., Linton, O., Lu, Z., 2018.  Semiparametric ultra-high dimensional model averaging of
nonlinear dynamic time series. Journal of the American Statistical Association, in press.

\vspace{1mm}
\noindent Cheng, X., Hansen, B., 2015. Forecasting with factor-augmented regression: a frequentist model averaging
approach. Journal of Econometrics, 186, 280--293.
\vspace{1mm}

\vspace{1mm}
\noindent Claeskens, G., Hjort, N. L., 2008. Model selection and model averaging. Cambridge University Press.
\vspace{1mm}

\vspace{1mm}
\noindent  Fan, J., Peng, H., 2004. On non-concave penalized likelihood with diverging number of
parameters. The Annals of Statistics, 32, 928--961.
\vspace{1mm}

\vspace{1mm}
\noindent Guo, Y., Zhang, X., Wang, S., Zou, G., 2016. Model averaging based on leave-subject-out cross-validation. Journal of Econometrics, 192, 139--151.
\vspace{1mm}

\vspace{1mm}
\noindent Hansen, B. E., 2007. Least squares model averaging. Econometrica, 75, 1175--1189.
\vspace{1mm}

\vspace{1mm}
\noindent Hansen, B. E., 2008. Least squares forecast averaging. Journal of Econometrics, 146, 342--350.
\vspace{1mm}

\vspace{1mm}
\noindent Hansen, B. E., Racine, J. S., 2012. Jackknife model averaging. Journal of Econometrics, 167, 38--46.
\vspace{1mm}

\vspace{1mm}
\noindent Horowitz, J.L., 1998. Bootstrap methods for median regression models. Econometrica, 66, 1327--1351.
\vspace{1mm}

\vspace{1mm}
\noindent Lian, H., 2012. A note on the consistency of Schwarz¡¯s criterion in linear quantile regression with the SCAD penalty. Statistics and Probability Letters, 82, 1224--1228.
\vspace{1mm}

\vspace{1mm}
\noindent Hjort, N., Claeskens, G., 2003. Frequentist model average estimators. Journal of the American Statistical Association, 98, 879--899.
\vspace{1mm}

\vspace{1mm}
\noindent Hoeting, J. A., Madigan, D., Raftery, A. E., Volinsky, C. T., 1999. Bayesian model averaging: a tutorial. Statistical Science, 14, 382--417.
\vspace{1mm}

\vspace{1mm}
\noindent Horowitz, J. L., Lee, S., 2005. Nonparametric estimation of an additive quantile regression model. Journal of the American Statistical Association, 100, 1238--1249.
\vspace{1mm}

\vspace{1mm}
\noindent Huang, T. and Li, J. 2018. Semiparametric model average prediction in panel data analysis.
Journal of Nonparametric Statistics, 30, 125--144.
\vspace{1mm}

\vspace{1mm}
\noindent Jiang, J., Zhao, Q., Hui, Y. V., 2001. Robust modelling of ARCH models. Journal of Forecasting,
20, 111--133.
\vspace{1mm}

\vspace{1mm}
\noindent Johnson, R.W., 1996. Fitting percentage of body fat to simple body measurements. Journal of Statistics Education, 4, 1--8.
\vspace{1mm}

\vspace{1mm}
\noindent Kai, B., Li, R., Zou, H., 2011. New efficient estimation and variable selection methods for semiparametric varying-coefficient partially linear models. The Annals of Statistics, 39, 305--332.
\vspace{1mm}

\vspace{1mm}
\noindent Knight, K., 1998. Limiting distributions for $L_1$ regression estimators under general conditions. The Annals of Statistics, 26,
755--770.
\vspace{1mm}

\vspace{1mm}
\noindent Li, D., Linton, O., Lu, Z., 2015. A flexible semiparametric forecasting model for time series. Journal of
Econometrics, 187, 345--357.
\vspace{1mm}

\vspace{1mm}
\noindent Li, J., Xia, X., Wong, W. K., Nott, D., 2018. Varying-coefficient semiparametric model averaging prediction. Biometrics, in press.
\vspace{1mm}

\vspace{1mm}
\noindent Liang, H., Zou, G., Wan, A. T. K., Zhang, X., 2011. Optimal weight choice for frequentist model average estimators. Journal of the American Statistical Association, 106, 1053--1066.
\vspace{1mm}

\vspace{1mm}
\noindent Liu, C.-A., 2015. Distribution theory of the least squares averaging estimator. Journal of Econometrics, 186, 142--159.
\vspace{1mm}

\vspace{1mm}
\noindent Lv, J., Yang, H., Guo, C.,  Smoothing combined generalized estimating equations quantile partially linear additive models longitudinal data. Computational Statistics, 31, 1203--1234.
\vspace{1mm}

\vspace{1mm}
\noindent Parzen, E., 1962. On estimation of a probability density function and model. The Annals of Statistics, 33, 1065--1076.
\vspace{1mm}

\vspace{1mm}
\noindent Pollard, D., 1991. Asymptotics for least absolute deviation regression estimators. Econometr. Theory, 7, 186--199.
\vspace{1mm}

\vspace{1mm}
\noindent Raftery, A., Madigan, D., Hoeting, J., 1997. Bayesian model averaging for linear regression models. Journal of the American Statistical Association, 92, 179--191.
\vspace{1mm}

\vspace{1mm}
\noindent Tang, Y., Song, X., Wang, H. J., Zhu, Z., 2013. Variable selection in high-dimensional quantile varying coefficient models.
Journal of Multivariate Analysis, 122, 115--132.
\vspace{1mm}

\vspace{1mm}
\noindent Wan, T. K., Zhang, X., Zou, G., 2010. Least squares model averaging by Mallows criterion. Journal of Econometrics, 156, 277--283.
\vspace{1mm}

\vspace{1mm}
\noindent Zhang, X., Wan, A. T. K., Zhou, S. Z., 2012. Focused information criteria, model selection,
and model averaging in a tobit model with a nonzero threshold. Journal of Business $\&$
Economic Statistics, 30, 132--142.
\vspace{1mm}

\vspace{1mm}
\noindent Zhang, X., Wan, A. T. K., Zou, G., 2013. Model averaging by Jackknife criterion in models with dependent data. Journal of Econometrics, 174, 82--94.
\vspace{1mm}

\vspace{1mm}
\noindent Zhang, X., Zou, G., Liang, H., 2014. Model averaging and weight choice in linear mixed-effects models. Biometrika, 101, 205--218.
\vspace{1mm}

\vspace{1mm}
\noindent Zhang, X., Yu, D., Zou, G., Liang, H., 2016. Optimal model averaging estimation for generalized linear models
and generalized linear mixed-effects models. Journal of American Statistical Association, In Press.
\vspace{1mm}

\bibliographystyle{model5-names}
\bibliography{<your-bib-database>}

%% not want to use model5-names.bst.



%% \bibitem must have one of the following forms:
%%   \bibitem[Jones et al.(1990)]{key}...
%%   \bibitem[Jones et al.(1990)Jones, Baker, and Williams]{key}...
%%   \bibitem[Jones et al., 1990]{key}...
%%   \bibitem[\protect\citeauthoryear{Jones, Baker, and Williams}{Jones
%%       et al.}{1990}]{key}...
%%   \bibitem[\protect\citeauthoryear{Jones et al.}{1990}]{key}...
%%   \bibitem[\protect\astroncite{Jones et al.}{1990}]{key}...
%%   \bibitem[\protect\citename{Jones et al., }1990]{key}...
%%   \harvarditem[Jones et al.]{Jones, Baker, and Williams}{1990}{key}...
%%


% \end{thebibliography}

\end{document}